\documentclass[12pt]{article}
\usepackage{amsmath,amsfonts,amssymb,amsthm,enumerate,mathrsfs}

\newcommand{\R}{{\mathbb R}}
\newcommand{\Z}{{\mathbb Z}}

\newcommand{\C}{{\mathbb C}}

\newcommand{\1}{{\bf 1}}
\newcommand{\F}{{\mathbb F}}

\newcommand{\wt}{{\rm wt}}
\newcommand{\Aut}{{\rm Aut}}

\newcommand{\eqa}{\begin{eqnarray}}
\newcommand{\eeqa}{\end{eqnarray}}
\newcommand{\eqn}{\begin{eqnarray*}}
\newcommand{\eeqn}{\end{eqnarray*}}

\newcommand{\Hom}{{\rm Hom}}

\newcommand{\Sym}{{\rm Sym}}

\newtheorem{dfn}{Definition}[section]
\newtheorem{pro}[dfn]{Proposition}
\newtheorem{thm}[dfn]{Theorem}

\newtheorem{lem}[dfn]{Lemma}
\newtheorem{cor}[dfn]{Corollary}
\newtheorem{rem}[dfn]{Remark}
\newtheorem{note}[dfn]{Note}

\newcommand{\Ker}{{\rm Ker\ }}

\newcommand{\NO}{\,{\raise0.25em\hbox{$\mathop{\hphantom{\cdot}}\limits^{_{\circ}}_{^{\circ}}$}}\,}

\newcommand{\qe}{\qed\vskip2ex}
\makeatletter
\@addtoreset{equation}{section}

\makeatother

\def\bl{\begin{lem}}
\def\el{\end{lem}}
\def\bt{\begin{thm}}
\def\et{\end{thm}}
\def\bp{\begin{pro}}
\def\ep{\end{pro}}
\def\br{\begin{rem}}
\def\er{\end{rem}}
\def\bc{\begin{cor}}
\def\ec{\end{cor}}
\def\bd{\begin{dfn}\rm}
\def\ed{\end{dfn}}
\def\bn{\begin{note}\rm}
\def\en{\end{note}}
\def\proof{{\it Proof.}}

\title{\begin{flushright}
\end{flushright}\Large An $E_8$-approach to the moonshine vertex operator algebra}
\author{Hiroki SHIMAKURA\footnote{The author was partially supported by Grants-in-Aid for Scientific Research (No. 20549004) and Excellent Young Researcher Overseas Visit Program, Japan Society for the Promotion of Science.}}
\date{\small{\it Department of Mathematics,\\
Aichi University of Education,\\
1 Hirosawa, Igaya-cho, Kariya-city, Aichi, 448-8542 Japan}\\
{\rm e-mail: shima@auecc.aichi-edu.ac.jp}\\
\vspace{0.5cm}
{\rm 2000} {\it Mathematics Subject Classification}. Primary  17B69; Secondary 20D08.
}

\begin{document}
\maketitle

\begin{abstract}
In this article, we study the moonshine vertex operator algebra starting with the tensor product of three copies of the vertex operator algebra $V_{\sqrt2E_8}^+$, and describe it by the quadratic space over $\F_2$ associated to $V_{\sqrt2E_8}^+$.
Using quadratic spaces and orthogonal groups, we show the transitivity of the automorphism group of the moonshine vertex operator algebra on the set of all full vertex operator subalgebras isomorphic to the tensor product of three copies of $V_{\sqrt2E_8}^+$, and determine the stabilizer of such a vertex operator subalgebra.
Our approach is a vertex operator algebra analogue of ``An $E_8$-approach to the Leech lattice and the Conway group" by Lepowsky and Meurman.
Moreover, we find new analogies among the moonshine vertex operator algebra, the Leech lattice and the extended binary Golay code.
\end{abstract}

\section*{Introduction}

The moonshine vertex operator algebra (VOA) is one of the most important VOAs.
A reason is that its automorphism group is isomorphic to the Monster, the largest sporadic finite simple group.
Hence the Monster can be studied as symmetries of the moonshine VOA.

The moonshine VOA was constructed by Frenkel-Lepowsky-Meurman \cite{FLM} from the Leech lattice as a VOA analogue of the construction of the  Leech lattice from the extended binary Golay code.
Hence, it is natural to regard VOAs as analogues of binary codes and lattices.
In fact, by this approach, the upper bound of the minimum conformal weight of a holomorphic VOA was given in \cite{Ho95}, and the notion of a conformal design based on a VOA was introduced in \cite{Ho08}.
For other example, some isomorphism problems of VOAs were solved in \cite{Sh5}.

Lepowsky and Meurman studied the Leech lattice in \cite{LM} starting with the orthogonal direct sum of three copies of the $E_8$ root lattice, and described it by the quadratic space over $\F_2$ associated to the $E_8$ root lattice.
This is a lattice analogue of Turyn's construction of the extended binary Golay code (cf.\ \cite{MS}).
These are useful for understanding the automorphism groups of the Golay code and the Leech lattice, the Mathieu group and the Conway group (\cite{CS,LM,Gr2}).
Hence it is really important for the study of the Monster to find a VOA analogue of these descriptions.

\medskip

In this article, we study the moonshine VOA starting with the tensor product of three copies of the VOA $V_{\sqrt2E_8}^+$, and describe it as a simple current extension of $(V_{\sqrt2E_8}^+)^{\otimes3}$ by using a certain maximal totally singular subspace of the quadratic space over $\F_2$ associated to $(V_{\sqrt2E_8}^+)^{\otimes3}$.
This is a VOA analogue of Lepowsky and Meurman's description of the Leech lattice.
Note that the moonshine VOA was already constructed by Miyamoto \cite{Mi04} as a simple current extension of $(V_{\sqrt2E_8}^+)^{\otimes3}$ without quadratic spaces.

An advantage of our description is that we can effectively apply results of quadratic spaces and orthogonal groups to VOAs.
The transitivity of the automorphism group of the moonshine VOA on the set of all full subVOAs  isomorphic to $(V_{\sqrt2E_8}^+)^{\otimes3}$ is proved by the uniqueness of certain maximal totally singular subspaces of the quadratic space, and a subgroup of the automorphism group of the moonshine VOA of shape $2^{15}.(2^{20}:({\rm SL}_5(2)\times\Sym_3))$ is described as a lift of a certain subgroup of the orthogonal group.
According to \cite{ATLAS}, a subgroup with this shape is a maximal $2$-local in the Monster.

Another advantage is that we can find some aspects of the moonshine VOA and the Monster analogous to the Leech lattice and the Conway group, and to the extended binary Golay code and the Mathieu group respectively.
For example, the transitivity in the previous paragraph is a VOA analogue of that of the Mathieu group on the set of all trios of the extended binary Golay code (cf.\ \cite{CS}) and that of the Conway group on the set of all sublattices of the Leech lattice isomorphic to $(\sqrt2E_8)^{\oplus3}$ (\cite{Gr2}).
Moreover, the similarity of the following shapes of the stabilizers of a trio in the Mathieu group, of $(\sqrt2E_8)^{\oplus3}$ in the Conway group, and of $(V_{\sqrt2E_8}^+)^{\otimes3}$ in the Monster is explained by the associated orthogonal groups and quadratic spaces:
$$2^6:({\rm SL}_3(2)\times\Sym_3),\ 2^3.(2^{12}:({\rm SL}_4(2)\times\Sym_3)),\ 2^{15}.(2^{20}:({\rm SL}_5(2)\times\Sym_3)).$$
As another example, we find the following equation for the dimension of the weight $2$ subspace of the moonshine VOA (Proposition \ref{PEqV}):
$$ (3\times156)+(3\times (2^5-1))\times 1^0\times8^2+(3\times (2^5-1)\times 2^8)\times 1^2\times 8^1=196884,$$
where $156$ is the dimension of the weight $2$ subspace of the VOA $V_{\sqrt2E_8}^+$.
This is a VOA analogue of the following equations for the number of octads in the extended binary Golay code and that of vectors of norm $4$ in the Leech lattice:
$$\begin{array}{ccccccc}
(3\times 1)&+&(3\times (2^3-1))\times 1^0\times2^2&+&(3\times (2^3-1)\times 2^4)\times 1^2\times 2^1&=&759,\\
(3\times240)&+&(3\times (2^4-1))\times 2^0\times16^2&+&(3\times (2^4-1)\times 2^6)\times 2^2\times 16^1&=&196560.
\end{array} $$

\medskip

Let us explain the main results.
First, we discuss quadratic spaces over $\F_2$, which is crucial for this article.
Let $(R,q)$ be a non-singular quadratic space of plus type over $\F_2$ and let $w:R\to\{0,1,2\}$ be the map defined by setting $$w(a)=\begin{cases}0\quad {\rm if}\ v=0,\\ 1\quad {\rm if}\ q(a)=1,\\ 2\quad {\rm if}\ q(a)=0\ {\rm and}\ v\neq0.\end{cases}$$
Let $(R^3,q^3)$ be the orthogonal direct sum  of three copies of $(R,q)$ and let $\mathcal{S}$ be a maximal totally singular subspace of $(R^3,q^3)$ satisfying $$ w^3(v)\ge 4,\quad \forall v\in \mathcal{S}\setminus\{0\},$$
where $w^3(a_1,a_2,a_2)=\sum_{i=1}^3 w(a_i)$.
Such a subspace was constructed in \cite{LM} from complementary maximal totally singular subspaces of $(R,q)$.
We show in Theorem \ref{TCh} that $\mathcal{S}$ is unique up to $\Aut(R^3,w^3)$, the subgroup of the orthogonal group $O(R^3,q^3)$ preserving $w^3$, and determine in Proposition \ref{PStab} the stabilizer of $\mathcal{S}$ in $\Aut(R^3,w^3)$.

Next, we describe the moonshine VOA by the quadratic space.
Let $V=V_{\sqrt2E_8}^+$ and let $R(V)$ denote the set of all isomorphism classes of irreducible $V$-modules.
By \cite{AD,ADL}, $R(V)$ has a $10$-dimensional vector space structure over $\F_2$ under the fusion rules.
By \cite{Sh2}, $R(V)$ has a non-singular quadratic form $q_V$ of plus type, and hence $(R(V),q_V)$ is a $10$-dimensional non-singular quadratic space of plus type over $\F_2$.
Let $\mathcal{S}$ be a maximal totally singular subspace of $(R(V)^3,q_V^3)$.
We identify $R(V)^3$ with $R(V^{3})$ (\cite{FHL,DMZ}), where $V^3$ is the tensor product of three copies of $V$.
Then by a theory of simple current extensions (\cite{Hu,Mi04,LY2}), the $V^{3}$-module  $\mathfrak{V}(\mathcal{S})=\oplus_{[M]\in\mathcal{S}}M$ has a holomorphic VOA structure, where $[M]$ is the isomorphism class of $M$.
In this case, for $[M]\in R(V^3)$, the lowest weight of $M$ is $w^3([M])/2$.
Assume that $\mathcal{S}$ satisfies the condition on $w^3$ in the previous paragraph.
Then the the weight $1$ subspace of $\mathfrak{V}(\mathcal{S})$ is trivial.
Since $V$ is framed, so is $\mathfrak{V}(\mathcal{S})$.
Hence by \cite{LY}, $\mathfrak{V}(\mathcal{S})$ is isomorphic to the moonshine VOA (Theorem \ref{MTVm}).

Finally, we study the moonshine VOA and the automorphism group by using the orthogonal group and the quadratic space.
By \cite{Sh2}, the stabilizer of $\mathcal{S}$ in $\Aut(V^3)(\cong \Aut(R(V)^3,w^3))$ lifts to a subgroup of $\Aut(\mathfrak{V}(\mathcal{S}))$, and its shape is $2^{15}.(2^{20}:({\rm SL}_5(2)\times\Sym_3))$.
The uniqueness of $\mathcal{S}$ up to $\Aut(R(V)^3,w^3)$ implies the transitivity of $\Aut(\mathfrak{V}(\mathcal{S}))$ on the set of all full subVOAs of $\mathfrak{V}(\mathcal{S})$ isomorphic to $(V_{\sqrt2E_8}^+)^{3}$ (Theorem \ref{MT}).
Moreover, we describe in Proposition \ref{PEqV} the dimension of the weight $2$ subspace in terms of quadratic spaces.

\medskip

In order to compare VOAs to binary codes and lattices, we study the extended binary Golay code and the Leech lattice by similar approaches in Section 3, which were already discussed in \cite{MS,CS,LM,Gr2,Gr}.

\paragraph{Notations}
\begin{small}
\begin{center}
\begin{tabular}{ll}
$\langle\ , \rangle$& The symplectic form on $R$ or on $R^k$.\\
$(\ ,\ )$& A positive definite bilinear form on $\R^n$ or the standard inner product on $\F_2^n$.\\
$\times$& The fusion rules for a VOA.\\
$2^{n}$& An elementary abelian $2$-group of order $2^n$.\\
$A.B$& A group extension with normal subgroup $A$ and quotient $B$.\\
$A:B$& A split extension with normal subgroup $A$ and quotient $B$.\\
$C^k$& The direct sum of $k$ copies of a binary code $C$.\\
$\mathfrak{C}(\mathcal{S})$& The binary code associated to $\mathcal{S}\subset R(C)^k$.\\ 
$E_8$& The $E_8$ root lattice.\\
$\Phi, \Psi$& Maximal totally singular subspaces of $R$, $R(C)$, $R(L)$ or $R(V)$.\\
$\Phi_{(ij)}$& $\Phi_{(ij)}=\{(a_1,a_2,\dots,a_k)\in R^k\mid a_i=a_j\in\Phi,\ a_l=0\ {\rm if}\ l\neq i,j\}$.\\
$\Psi_{(12\dots k)}$& $\Psi_{(12\dots k)}=\{(b,\dots,b)\in R^k\mid b\in\Psi\}.$\\
$(\Phi\cap\Psi)_{(1)}$&$(\Phi\cap\Psi)_{(1)}=\{(a,0,\dots,0)\in R^k\mid a\in\Phi\cap\Psi\}.$\\
$g\circ M$& The conjugate of a module $M$ for a VOA by an automorphism $g$.\\
$g\circ [M]$& The isomorphism class of $g\circ M$.\\
$L^k$& The orthogonal direct sum of $k$ copies of a lattice $L$.\\
$\mathfrak{L}(\mathcal{S})$& The lattice associated to $\mathcal{S}\subset R(L)^k$.\\ 

$[M]$& The isomorphism class of a module $M$ for a VOA.\\
$O(R,q)$& The orthogonal group of a quadratic space $(R,q)$.\\
$O_2(G)$& The maximal normal $2$-subgroup of a group $G$.\\
$\rho_i$& The $i$-th coordinate projection from $R^k$ to $R$.\\
$q$& A quadratic form on $R$ of plus type\\
$q_C$& The quadratic form on $R(C)$ defined by $q_C(x)=\wt(x)/2\pmod 2$.\\
\end{tabular}
\end{center}
\begin{center}
\begin{tabular}{ll}
$q_L$& The quadratic form on $R(L)$ defined by $q_L(x)=( x,x)\pmod 2$.\\
$q_V$& The quadratic form on $R(V)$ defined by\\
& $q_V([M])=0$, $1$ if $M$ is $\Z$-graded, $(\Z+1/2)$-graded respectively.\\
${\rm SL}_n(2)$& The special linear group of degree $n$ over $\F_2$.\\
$\Sym_n$& The symmetric group of degree $n$.\\
$\mathcal{S}$& A maximal totally singular subspace of $R^k$.\\
$\mathcal{S}(\Phi,\Psi;k)$& The maximal totally singular subspace of $R^k$\\ & spanned by $\Phi_{(1i)}$ $(2\le i\le k)$, $\Psi_{(12\dots k)}$ and $(\Phi\cap\Psi)_{(1)}$.\\
$(R,q)$& A non-singular quadratic space of plus type with $q$ over $\F_2$.\\
$(R^k,q^k)$& The orthogonal direct sum of $k$ copies of $(R,q)$.\\
$R(C)$& $C^\perp/C$, where $C^\perp$ is the dual code of $C$.\\
$R(L)$& $L^*/L$, where $L^*$ is the dual lattice of $L$.\\
$R(V)$& The set of all isomorphism classes of irreducible modules for a VOA $V$.\\
$V$& A simple rational $C_2$-cofinite self-dual VOA of CFT type, or $V=V_{\sqrt2E_8}^+$.\\
$V^k$& The tensor product of $k$ copies of a VOA $V$.\\
$\mathfrak{V}(\mathcal{S})$& The VOA associated to $\mathcal{S}\subset R(V)^k$.\\ 
$w^k$& The map from $R^k$ to $\{0,1,\dots,2k\}$ (see Section 2.1).\\
$\wt(x)$& The (Hamming) weight of $x=(x_i)\in\F_2^n$ defined by $\wt(x)=|\{i\mid x_i\neq0\}|$.\\
\end{tabular}
\end{center}
\end{small}

\section{Preliminary}
In this section, we recall or give some definitions and facts necessary in this article.

\subsection{Quadratic spaces and orthogonal groups}
Let us recall fundamental facts on quadratic spaces over $\F_2$ and orthogonal groups (cf.\ \cite{Ch}).

Let $R$ be a $2m$-dimensional vector space over $\F_2$.
A form $\langle\ ,\ \rangle:R\times R\to \F_2$ is said to be {\it symplectic} if it is a symmetric bilinear form with $\langle a,a\rangle=0$ for all $a\in R$.
A map $q:R\to\F_2$ is called a {\it quadratic} form on $R$ if the associated form defined by $\langle a,b\rangle=q(a+b)+q(a)+q(b)$, $a,b\in R$, is symplectic.
A quadratic form is said to be {\it non-singular} if the associated symplectic form is non-singular, that is, $R^\perp=\{a\in R\mid \langle a,R\rangle=0\}=0$.
The pair $(R,q)$ consisting of a vector space $R$ over $\F_2$ and a quadratic form $q$ on it is called a {\it quadratic space}, and it is said to be {\it non-singular} if $q$ is non-singular.
A vector $a\in R$ is said to be {\it singular} (resp. {\it non-singular}) if $q(a)=0$ (resp. $q(a)=1$).
A subspace $\Phi$ of $R$ is said to be {\it totally singular} if any vector in $\Phi$ is singular.
A non-singular quadratic form $q$ is said to be {\it of plus type} if the dimension of a maximal totally singular subspace of $(R,q)$ is $m$.
Let $O(R,q)$ denote the orthogonal group of $(R,q)$, the subgroup of ${\rm GL}(R)$ preserving $q$.
The following lemma is well-known.

\bl\label{LO} {\rm (cf.\ \cite{Ch})} Let $(R,q)$ be a non-singular $2m$-dimensional quadratic space of plus type over $\F_2$.
Let $\Phi$ be a maximal totally singular subspace of $R$ and let $H$ be the stabilizer of $\Phi$ in $O(R,q)$.
\begin{enumerate}
\item The subspace of $R$ orthogonal to $\Phi$ is equal to $\Phi$, that is, $\Phi^\perp=\Phi$.
\item The orthogonal group $O(R,q)$ is transitive on the set of all maximal totally singular subspaces of $R$.
\item The maximal normal $2$-subgroup $O_2(H)$ of $H$ acts trivially on both $\Phi$ and $R/\Phi$, and $H/O_2(H)$ acts on $\Phi$ as ${\rm SL}(\Phi)(\cong{\rm SL}_m(2))$.
Moreover, $H$ has the shape $2^{\binom{m}{2}}:{\rm SL}_m(2)$.
\item There exists a maximal totally singular subspace $\Psi$ of $R$ such that $\Phi\cap \Psi=0$.
Moreover, the stabilizer of $\Psi$ in $H$ is isomorphic to ${\rm SL}_m(2)$. 
\item All pairs of complementary maximal totally singular subspaces of $R$ are conjugate under $O(R,q)$. 
\item $O(R,q)$ is transitive on the set of all non-zero singular vectors in $R$.
\end{enumerate}
\el
\subsection{Vertex operator algebras}
In this subsection, we recall the definitions of vertex operator algebras (VOAs) and modules from \cite{Bo,FLM,FHL}.
Throughout this article, all VOAs are defined over the field $\C$ of complex numbers unless otherwise stated.

A vertex operator algebra (VOA) $(V,Y,\1,\omega)$ is a $\Z_{\ge0}$-graded
 vector space $V=\oplus_{m\in\Z_{\ge0}}V_m$ equipped with a linear map

$$Y(a,z)=\sum_{i\in\Z}a_iz^{-i-1}\in ({\rm End}(V))[[z,z^{-1}]],\quad a\in V,$$
the {\it vacuum vector} $\1$ and the {\it Virasoro element} $\omega$ satisfying a number of conditions (\cite{Bo,FLM}).
We often denote it by $V$ or $(V,Y)$ simply.

Two VOAs $(V,Y,\1,\omega)$ and $(V^\prime,Y',\1',\omega')$ are said to be {\it isomorphic} if there exists a linear isomorphism $g$ from $V$ to $V^\prime$ satisfying $$ g\omega=\omega'\quad {\rm and}\quad gY(v,z)=Y'(gv,z)g\quad {\rm for}\ \forall v\in V.$$ 
When two VOAs are the same, such a linear isomorphism is called an {\it automorphism}.
The group of all automorphisms of $V$ is called the {\it automorphism group} of $V$ and is denoted by $\Aut(V)$.
We mean by a {\it subVOA} (or a {\it vertex operator subalgebra}) a graded subspace of $V$ which has a structure of a VOA such that the operations and its grading agree with the restriction of those of $V$ and that they share the vacuum vector.
When they share also the Virasoro element, we will call it a {\it full subVOA}.

An (ordinary) module $(M,Y_M)$ for a VOA $V$ is a $\C$-graded vector space $M=\oplus_{m\in\C} M_{m}$ equipped with a linear map
$$Y_M(a,z)=\sum_{i\in\Z}a_iz^{-i-1}\in ({\rm End}(M))[[z,z^{-1}]],\quad a\in V$$
satisfying a number of conditions (\cite{FHL}).
We often denote it by $M$ and its isomorphism class by $[M]$.
The {\it weight} of a homogeneous vector $v\in M_k$ is $k$.
If $M$ is irreducible then there exists the {\it lowest weight} $h\in\C$ of $M$ such that $M=\oplus_{m\in\Z_{\ge0}}M_{h+m}$ and $M_h\neq0$.

A VOA $V$ is said to be {\it of CFT type} if $V_0=\C\1$, is said to be {\it rational} if any module is completely reducible, and is said to be {\it $C_2$-cofinite} if $V/{\rm Span}_\C\{ a_{-2}b\mid a,b\in V\}$ is finite-dimensional.
A VOA is said to be {\it holomorphic} if it is the only irreducible module up to isomorphism.
A module $M$ is said to be {\it self-dual} if its contragredient module (cf.\ \cite{FHL}) is isomorphic to itself.
Let $R(V)$ denote the set of all isomorphism classes of irreducible $V$-modules.
Note that if $V$ is rational then $|R(V)|<\infty$.

Let $M_a,M_b,M_c$ be modules for a simple rational $C_2$-cofinite VOA $V$.
An intertwining operator of type $M_a\times M_b\to M_c$ is a linear map $M_a\mapsto ({\rm Hom}(M_b,M_c))\{z\}$ satisfying a number of conditions (\cite{FHL}).
The {\it fusion rule} $N_{M_a,M_b}^{M_c}$ is the dimension of the space of all intertwining operators of type $M_a\times M_b\to M_c$.
Since $V$ is $C_2$-cofinite, the fusion rules are finite (\cite{ABD}).
By the definition of the fusion rules, $N_{M_a,M_b}^{M_c}=N_{M_a',M_b'}^{M_c'}$ if $M_x\cong M_x'$ as $V$-modules for all $x=a,b,c$.
Hence, the fusion rules are determined by the isomorphism classes of $V$-modules.
For irreducible modules $M_a$ and $M_b$, we use the following expression $$[M_a]\times [M_b]=\sum_{[M]\in R(V)} N_{M_a,M_b}^{M}[M],$$
which is also called the fusion rule.

Let $M=(M,Y_M)$ be a module for a VOA $V$.
For $g\in\Aut(V)$, let $g\circ M=(M,Y_{g\circ M})$ denote the $V$-module defined by $Y_{g\circ M}(v,z)=Y_M(g^{-1}v,z)$.
\bl\label{LAct0} Let $M$ and $M'$ be $V$-modules and let $g\in\Aut(V)$.
\begin{enumerate}
\item If $M\cong M'$ as $V$-modules then $g\circ M\cong g\circ M'$ as $V$-modules.
\item If $M$ is $(h+\Z_{\ge0})$-graded then so is $g\circ M$.
\end{enumerate}
\el
\proof\ Let $f:M\to M'$ be an isomorphism of $V$-modules, that is, $fY_M(v,z)w=Y_{M'}(v,z)f(w)$ for all $v\in V$ and $w\in M$.
Replacing $v$ by $g^{-1}v$, we obtain $fY_{g\circ M}(v,z)w=Y_{g\circ M'}(v,z)f(w)$ for all $v\in V$.
Hence $g\circ M\cong g\circ M'$ as $V$-modules, and we obtain (1).

By the axioms of modules, $\omega_1$ acts by $m$ on $M_m$.
Hence (2) follows from $g^{-1}(\omega)=\omega$.\qe

If $M$ is irreducible then so is $g\circ M$.
By the lemma above, $\Aut(V)$ acts on $R(V)$.
For $\mathcal{T}=\{W^i\mid 1\le i\le s\}\subset R(V)$, let $g\circ\mathcal{T}$ denote $\{g\circ W^i\mid 1\le i\le s\}$. 

\bl\label{LAct} Let $V$ be a simple rational $C_2$-cofinite VOA and let $M_a,M_b,M_c$ be $V$-modules.
Let $g\in\Aut(V)$.
Then $N_{g\circ M_{a},g\circ M_{b}}^{g\circ M_{c}}=N_{M_{a},M_{b}}^{M_{c}}$.
In particular, the action of $\Aut(V)$ on $R(V)$ preserves the fusion rules, that is, for $[M_a],[M_b]\in R(V)$,
$$(g\circ [M_{a}])\times (g\circ [M_{b}])=\sum_{[M]\in R(V)}N_{M_{a},M_{b}}^{M}(g\circ [M]).$$
\el
\proof\ By the definition of $g\circ M_x$, for each $x\in\{a,b,c\}$, there exists a non-zero linear isomorphism $\nu_x:g\circ M_x\to M_x$ such that $\nu_xY_{g\circ M_{x}}(v,z)=Y_{M_{x}}(g^{-1}v,z)\nu_x$ for $v\in V$.
Let $\mathcal{Y}$ be a non-zero intertwining operator of type $M_a\times M_b\to M_c$.
Then one can verify that 
\begin{eqnarray*}
(g\circ \mathcal{Y})(v,z)=\nu_c^{-1}\mathcal{Y}(\nu_av,z)\nu_b\in(\Hom(g\circ M_b,g\circ M_c))\{z\},\ v\in g\circ M_a,
\end{eqnarray*}
is a non-zero intertwining operator of type $(g\circ M_a)\times (g\circ M_b)\to (g\circ M_c)$.
Hence $N_{M_{a},M_{b}}^{M_c}\le N_{g\circ M_{a},g\circ M_{b}}^{g\circ M_c}$.
Moreover, $N_{g\circ M_{a},g\circ M_{b}}^{g\circ M}\le N_{M_{a},M_{b}}^{M_c}$ follows from $g^{-1}\circ (g\circ M)\cong M$.
Thus $N_{M_{a},M_{b}}^{M_c}= N_{g\circ M_{a},g\circ M_{b}}^{g\circ M_c}$.
\qe

The theory of tensor products of VOAs was established in \cite{FHL}.
For a positive integer $k$, let $V^k$ denote the tensor product of $k$ copies of $V$.
Later, we use the following lemma.
\bl\label{LemFHL} {\rm (\cite[Section 4.7]{FHL}, \cite{DMZ})} Let $V$ be a simple rational $C_2$-cofinite VOA of CFT type.
Then $$R({V^k})=\{\otimes_{i=1}^k W_i\mid W_i\in R(V)\},$$ and for $\otimes_{i=1}^k W_{i,a}, \otimes _{i=1}^kW_{i,b},\otimes _{i=1}^kW_{i,c}\in R(V^k)$, the following fusion rule holds: 
$$\otimes_{i=1}^k W_{i,a}\times \otimes _{i=1}^kW_{i,b}=\otimes_{i=1}^k (W_{i,a}\times W_{i,b}).$$
\el

\subsection{Simple current extensions of VOAs}
In this subsection, we recall the notion of simple current modules and simple current extensions.

Let $V(0)$ be a simple VOA.
An irreducible $V(0)$-module $M_a$ is called a {\it simple current module} if for any irreducible $V(0)$-module $M_b$, there exists a unique irreducible $V(0)$-module $M_c$ satisfying the fusion rule $[M_a]\times [M_b]=[M_c]$.
A simple VOA $V$ is called a {\it simple current extension} of $V(0)$ graded by a finite abelian group $E$ if $V$ is the direct sum of non-isomorphic irreducible simple current $V(0)$-modules $\{V(\alpha)\mid \alpha\in E\}$ and the fusion rule $[V(\alpha)]\times [V(\beta)]=[V(\alpha+\beta)]$ holds for all $\alpha,\beta\in E$.
The uniqueness theorem for simple current extensions was proved in \cite{DM}.

\bp\label{PSC} {\rm \cite[Proposition 4.3]{DM}} Let $(V,Y)$ be a simple current extension of a simple VOA $V(0)$.
Then the VOA structure of $V$ as a simple current extension of $V(0)$ is uniquely determined by the $V(0)$-module structure of $V$, that is, if $({V}',{Y}')$ has a VOA structure with $V={V}'$ as vector spaces and $Y(v,z)={Y}'(v,z)$ for all $v\in V(0)$ then the VOAs $(V,Y)$ and $({V}',{Y}')$ are isomorphic.
\ep

The conjugates of a simple current extension by automorphisms were studied in \cite{SY}.

\bl\label{LSY} {\rm \cite[Lemma 3.14]{SY}} Let $V=\oplus_{\alpha\in E}V(\alpha)$ be a simple current extension of a simple VOA $V(0)$ graded by a finite abelian group $E$.
Let $g$ be an automorphism of $V(0)$.
Then the $V(0)$-module $\oplus_{\alpha\in E}\ g\circ V(\alpha)$ has a VOA structure isomorphic to $V$.
\el

Let $V=\oplus_{\alpha\in E}V(\alpha)$ be a simple current extension of a simple VOA $V(0)$ graded by a finite abelian group $E$.
Then the dual $E^*$of $E$ acts on $V$ as an automorphism group: $\chi\in E^*$ acts on $V(\alpha)$ by the scalar multiplication $\chi(\alpha)$ for each $\alpha\in E$.
The following restriction homomorphism was studied in \cite{Sh2}:
\begin{eqnarray}
 \eta:N_{\Aut(V)}(E^*)\to\biggl\{ g\in\Aut(V(0))\biggl|\ g\circ \{[V(\alpha)]\mid \alpha\in E\}=\{[V(\alpha)]\mid \alpha\in E\}\biggr\}.\label{eta}
\end{eqnarray}
\bp\label{PSh2} {\rm \cite[Theorem 3.3]{Sh2}} The map $\eta$ is surjective, and $\Ker\eta=E^*$.
\ep

Let us discuss the group structure on a subset of $R(V)$ under the fusion rules.

\bl\label{LGp}{\rm (cf.\ \cite[Corollary 1 (3)]{LY2})} Let $V$ be a simple rational $C_2$-cofinite VOA of CFT type and $\mathcal{T}$ a subset of $R(V)$.
Assume that a representative of any element in $\mathcal{T}$ is a self-dual simple current module and that $\mathcal{T}$ is closed under the fusion rules, that is, $W \times W'=W''\in \mathcal{T}$ for all $W,W'\in \mathcal{T}$.
Then $(\mathcal{T},\times)$ is an elementary abelian $2$-group.
\el
\proof\ By the assumption on $\mathcal{T}$, the fusion rule $\times$ is a binary operation on $\mathcal{T}$.
By the assumptions on $V$, the fusion rules are associative (\cite{Hu2}).
For any $[M]\in \mathcal{T}$, $[M]\times [M]=[V]$ since $[M]$ is self-dual.
Hence $[V]\in\mathcal{T}$, and clearly $[V]$ is the identity.
Thus $(\mathcal{T},\times)$ is an elementary abelian $2$-group.\qe

The existence theorem for a simple VOA structure on a direct sum of non-isomorphic self-dual simple current modules was established in \cite{Hu,LY2}.
For the definition of invariant bilinear forms, see \cite{FHL}.

\bp\label{PVOA} {\rm (\cite{Hu}, \cite[Proposition 3, Theorem 2]{LY2}) (cf.\ \cite{Mi04})} Let $V$ be a simple rational $C_2$-cofinite VOA of CFT type and let $\mathcal{T}$ be a set of the isomorphism classes of non-isomorphic self-dual $\Z$-graded simple current $V$-modules.
Assume that for any element $[M]\in\mathcal{T}$ the invariant bilinear form on $M$ is symmetric.
Then $V$-module $\oplus_{[M]\in \mathcal{T}}M$ has a simple VOA structure by extending its $V$-module structure if and only if $\mathcal{T}$ is closed under the fusion rules.
\ep

\subsection{VOA $V_{\sqrt2E_8}^+$}
Let $E_8$ denote the $E_8$ root lattice and set $\sqrt2E_8=\{\sqrt2v\mid v\in E_8\}$.
Let $V_{\sqrt2E_8}$ be the lattice VOA associated with $\sqrt2E_8$ (\cite{Bo,FLM}) and let $V_{\sqrt2E_8}^+$ be the subspace of $V_{\sqrt2E_8}$ fixed by the automorphism of $V_{\sqrt2E_8}$ of order $2$ which is a lift of $-1\in\Aut(\sqrt2E_8)$.
Then $V_{\sqrt2E_8}^+$ is a simple full subVOA of $V_{\sqrt2E_8}$ of CFT type, and its central charge is $8$.
Since $\sqrt2E_8$ has an orthogonal basis of norm $4$, $V_{\sqrt2E_8}^+$ is a framed VOA.
By \cite{DGH}, $V_{\sqrt2E_8}^+$ is rational and $C_2$-cofinite.
In this subsection, we review the properties of $V=V_{\sqrt2E_8}^+$ and the set $R(V)$ of all isomorphism classes of irreducible $V$-modules.

In \cite{AD}, $R(V)$ was determined, and $|R(V)|=2^{10}$.
Moreover, in \cite{ADL}, the fusion rules of $R(V)$ were completely calculated and the contragredient modules of irreducible $V$-modules were determined.
In particular, any irreducible $V$-module is a self-dual simple current module.
By Lemma \ref{LGp}, $R(V)$ has an elementary abelian $2$-group structure of order $2^{10}$ under the fusion rules.
We view $R(V)$ as a $10$-dimensional vector space over $\F_2$.
By Lemma \ref{LAct}, $\Aut(V)$ acts on $R(V)$ as a subgroup of ${\rm GL}(R(V))$.

Let $q_{V}:R(V)\to \F_2$ be the map defined by setting $q_{V}([M])=1$ if $M$ is $(1/2+\Z)$-graded, and $q_{V}([M])=0$ if $M$ is $\Z$-graded.
Note that for $[M]\in R(V)$, the lowest weight of $M$ is $w([M])/2$, where $w([M])=0,1$ and $2$ if $[M]=[V]$, $q([M])=1$, and $q([M])=0$ with $[M]\neq[V]$ respectively (cf.\ Section 2.1).
It was shown in \cite[Theorem 3.8]{Sh2} that $q_{V}$ is a non-singular quadratic form on $R(V)$ of plus type. 
By Lemma \ref{LAct0} (2), $\Aut(V)$ preserves $q_V$.
Hence we obtain a group homomorphism from $\Aut(V)$ to the orthogonal group $O(R(V),q_{V})$.
In fact, it is an isomorphism by \cite[Theorem 4.5]{Sh2} (cf.\ \cite{Gr1}).

Since $V_{\sqrt2E_8}=V\oplus V_{\sqrt2E_8}^-$ is a VOA, the invariant bilinear form on the irreducible $V$-module $V_{\sqrt2E_8}^-$ is symmetric (\cite[Proposition 5.3.6]{FHL}).
By Lemma \ref{LO} (6), for any $\Z$-graded irreducible $V$-module, the invariant bilinear form on it is also symmetric.
Since the lowest weight space of any $(\Z+1/2)$-graded irreducible $V$-module is one-dimensional, the invariant form must be symmetric.

\bl\label{POVE} Let $V=V_{\sqrt2E_8}^+$.
\begin{enumerate}[{\rm (1)}]
\item $V$ is simple, rational, $C_2$-cofinite, self-dual and of CFT type.
\item $(R(V),q_V)$ is a non-singular $10$-dimensional quadratic space of plus type over $\F_2$.
\item For $[M]\in R(V)$, the lowest weight of $M$ is $w([M])/2$.
\item $\Aut(V)\cong O(R(V),q_V)$.
\item For any irreducible $V$-module, the invariant bilinear form on it is symmetric.
\end{enumerate}
\el

\section{Direct sum of quadratic spaces}

Let $(R,q)$ be a non-singular $2m$-dimensional quadratic space of plus type over $\F_2$ and let $k$ be a positive integer.
Let $R^k$ denote the orthogonal direct sum of $k$ copies of $R$.
We use the standard expression $(a_1,a_2,\dots,a_k)$ for a vector in $R^k$.
Let $q^k:R^k\to\F_2$, $(a_1,a_2,\dots,a_k)\mapsto\sum_{i=1}^kq(a_i)$.
Then $(R^k,q^k)$ is a non-singular $2mk$-dimensional quadratic space of plus type over $\F_2$. 
Note that the symmetric group $\Sym_k$ of degree $k$ acts naturally on $R^k$ as the permutation group on the coordinates.
Clearly, $\Sym_k\subset O(R^k,q^k)$.

In this section, we introduce a map $w^k:R^k\to\{0,1,\dots, 2k\}$ and study maximal totally singular subspaces $\mathcal{S}$ of $(R^k,q^k)$ such that $w^k(v)\ge4$ for all $v\in\mathcal{S}\setminus\{0\}$.

\subsection{A construction of maximal totally singular subspaces}
In this subsection, we construct a maximal totally singular subspace of $R^k$ from a pair of maximal totally singular subspaces of $R$, which is a slight generalization of \cite{LM} (cf.\ \cite{Gr}), and describe its stabilizer in $O(R,q)\wr\Sym_k$. 

Let $w:R\to\{0,1,2\}$ be the map defined by 
$$w(a)=\begin{cases}0\quad {\rm if}\ a=0,\\ 1\quad {\rm if}\ q(a)=1,\\ 2\quad {\rm if}\ q(a)=0\ {\rm and}\ a\neq0,\end{cases}$$
and let $$w^k:R^k\to\{0,1,\dots, 2k\},\quad (a_1,a_2,\dots,a_k)\mapsto\sum_{i=1}^k w(a_i).$$
Since $w(a)\equiv q(a)\pmod 2$ for all $a\in R$, we have $w^k(v)\equiv q^k(v)\pmod 2$ for all $v\in R^k$.
Hence $w^k(v)\in2\Z$ if and only if $v$ is singular.
The following lemma is easy.

\bl\label{Lw2} Let $v$ be a vector in $(R^k,q^k)$.
\begin{enumerate}[{\rm (1)}]
\item $w^k(v)=1$ if and only if $v=\sigma(a,0,\dots,0)$ for some non-singular vector $a\in R$ and $\sigma\in\Sym_k$.
\item $w^k(v)=2$ if and only if $v=\sigma(a,0,\dots,0)$ or $v=\sigma(b,c,0,\dots,0)$ for some non-zero singular vector $a\in R$ or non-singular vectors $b,c\in R$ and $\sigma\in \Sym_k$.
\end{enumerate}
\el

Let $\Aut(R^k,w^k)$ denote the subgroup of ${\rm GL}(R^k)$ preserving $w^k$.
Then $\Aut(R^k,w^k)$ is a subgroup of $O(R^k,q^k)$ since $w^k(v)\equiv q^k(v)\pmod 2$ for $v\in R^k$.

\bp\label{PAutw} $\Aut(R^k,w^k)\cong O(R,q)\wr\Sym_k$.
\ep
\proof\ Let $G=\Aut(R^k,w^k)$.
Clearly, $O(R,q)\wr\Sym_k$ acts faithfully on $R^k$ and preserves $w^k$.
Hence $O(R,q)\wr\Sym_k\subset G$.

Set $R_{(i)}=\{(a_1,a_2,\dots,a_k)\in R^k\mid a_j=0\ {\rm if}\ j\neq i\}$.
By Lemma \ref{Lw2} (1), \begin{equation}\{v\in R^k\mid w^k(v)=1\}=\{v\in \cup_{i=1}^kR_{(i)}\mid q^k(v)=1\}.\label{Eq:wt1}\end{equation}
Hence $G$ preserves $\{v\in \cup_{i=1}^kR_{(i)}\mid q^k(v)=1\}$.
Moreover, by Lemma \ref{Lw2} (2) and (\ref{Eq:wt1}), $G$ also preserves $\{v\in \cup_{i=1}^kR_{(i)}\mid q^k(v)=0\}$.
Hence $G$ preserves $\cup_{i=1}^k R_{(i)}$.

Let $a,b\in R_{(i)}$ with $\langle a,b\rangle=1$ and let $g\in G$.
If $g(a)\in R_{(j)}$ then $g(b)\in R_{(j)}$ since $\langle R_{(p)},R_{(q)}\rangle=0$ if $p\neq q$.
Since $R$ is non-degenerate, for non-zero $a\in R$, ${\rm Span}_{\F_2}\{b\in R\mid \langle b,a\rangle=1\}=R$. 
Hence $G$ permutes $\{R_{(i)}\mid i=1,2,\dots,k\}$.
It follows from $O(R_{(i)},q^k)\cong O(R,q)$ that $G\subset O(R,q)\wr\Sym_k$.\qe

Now, let us construct maximal totally singular subspaces $\mathcal{S}$ of $R^k$ satisfying
\begin{equation}
w^k(v)\ge 4,\quad \forall v\in \mathcal{S}\setminus\{0\}\label{Cond1}
\end{equation}
as a natural generalization of \cite{LM} (cf.\ \cite{Gr}).
Let $\Phi,\Psi$ be maximal totally singular subspaces of $R$.
For $\{i,j\}\subset\{1,2,\dots,k\}$, set 
$$\Phi_{(ij)}=\{(a_1,a_2,\dots,a_k)\in R^k\mid a_i=a_j\in\Phi,\ a_l=0\ {\rm if}\ l\neq i,j\}$$
and $$\Psi_{(12\dots k)}=\{(b,\dots,b)\in R^k\mid b\in\Psi\}.$$
Set $(\Phi\cap\Psi)_{(1)}=\{(a,0,\dots,0)\in R^k\mid a\in\Phi\cap\Psi\}$ and $\mathcal{S}(\Phi,\Psi;k)={\rm Span}_{\F_2}\{\Phi_{(1i)},\Psi_{(12\dots k)},(\Phi\cap\Psi)_{(1)}\mid 2\le i\le k\}$.
Then 
\begin{equation}
\mathcal{S}(\Phi,\Psi;k)=\{(a_1+b,a_2+b,\dots,a_{k}+b)\in R^k\mid a_i\in\Phi, b\in\Psi, \sum_{i=1}^ka_i\in\Phi\cap\Psi\}.\label{EqT}
\end{equation}
\br If $k$ is odd then $\mathcal{S}(\Phi,\Psi;k)={\rm Span}_{\F_2}\{\Phi_{(1i)},\Psi_{(12\dots k)}\mid 2\le i\le k\}$.
\er

\bp\label{PLM} {\rm (cf.\ \cite[p487]{LM})} The subspace $\mathcal{S}(\Phi,\Psi;k)$ of $(R^k,q^k)$ is maximal totally singular.
Moreover if $k\ge3$ and $\Phi\cap\Psi=0$ then it satisfies (\ref{Cond1}). 
\ep
\proof\ By the definition, $\Phi_{(12)},\Phi_{(13)},\cdots,\Phi_{(1k)}$ and $\Psi_{(12\dots k)}$ are mutually perpendicular, and $\Phi_{(1i)}$ and $\Psi_{(12\dots k)}$ are totally singular.
Hence $\mathcal{S}(\Phi,\Psi;k)$ is totally singular.
Let us count the dimension of $\mathcal{S}(\Phi,\Psi;k)$.
Set $U=\{(a_1,a_2,\dots,a_k)\in R^k\mid a_i\in\Phi\cap\Psi\}$ and $d=\dim \Phi\cap\Psi$.
Then $U\subset \mathcal{S}(\Phi,\Psi;k)$ and $\dim U=kd$.
Set $T=\mathcal{S}(\Phi,\Psi;k)/U$.
Then $T=\oplus_{i=2}^k\bar{\Phi}_{(1i)}\oplus\bar{\Psi}_{(12\dots k)}$ as vector spaces, where $\bar{\Psi}=\Psi/(\Phi\cap\Psi)$ and $\bar{\Phi}=\Phi/(\Phi\cap\Psi)$.
Since this is a direct sum, $\dim T= (k-1)\dim\bar{\Phi}+\dim\bar{\Psi}$.
Hence $$\dim \mathcal{S}(\Phi,\Psi;k)=\dim T+\dim U=(k-1)(m-d)+(m-d)+kd=mk.$$
Thus $\mathcal{S}(\Phi,\Psi;k)$ is a maximal totally singular subspace of $(R^k,q^k)$.

Assume that $k\ge3$ and that $\Phi\cap\Psi=0$.
Let $v=(a_1+b,a_2+b,\dots,a_{k-1}+b,\sum_{i=1}^{k-1} a_i+b)\in\mathcal{S}(\Phi,\Psi;k)$ with  $w^k(v)\le2$, where $a_i\in\Phi$ and $b\in\Psi$.
Then at least $(k-2)$ entries are zero by Lemma \ref{Lw2} (2).
Say $a_i+b=0$, $1\le i\le k-2$, which implies $a_i=b=0$ by $\Psi\cap\Phi=0$ and $k\ge3$.
Hence $v=(0,\dots,0,a_{k-1},a_{k-1})$.
Since $a_{k-1}$ is singular and $w^k(v)\le 2$, we obtain $a_{k-1}=0$ by Lemma \ref{Lw2} (2).
Therefore $v=0$.\qe

Let us describe the stabilizer of $\mathcal{S}(\Phi,\Psi;k)$ in $\Aut(R^k,w^k)$.
Let $O(R,q)^k$ denote the normal subgroup of $O(R,q)\wr\Sym_k$ isomorphic to the direct product of $k$ copies of the orthogonal group $O(R,q)$.
For an element $g\in O(R,q)^{k}$, we use the expression $g=(g_1,g_2,\dots,g_k)$, where $g_i\in O(R,q)$.
For a subgroup $F$ of $O(R,q)$, we denote $F_{(ij)}=\{(g_1,g_2,\dots,g_k)\in O(R,q)^k\mid g_i=g_j\in F,\ g_l=1\ {\rm if}\ l\neq i,j\}$ and $F_{(12\dots k)}=\{(h,\dots,h)\in O(R,q)^k\mid h\in F\}$.

\bp\label{PStab} Let $G$ be the stabilizer of $\mathcal{S}(\Phi,\Psi;k)$ in $\Aut(R^k,w^k)$.
Let $H$ be the stabilizer of $\Phi$ in $O(R,q)$ and let $K$ be the stabilizer of $\Psi$ in $H$.
If $k\ge3$ and $\Phi\cap\Psi=0$ then $G$ is generated by $O_2(H)_{(1i)}(\cong 2^{\binom{m}{2}})$, $(2\le i\le k)$, $K_{(12\dots k)}(\cong{\rm SL}_m(2))$ and $\Sym_k$, and $G$ has the shape $2^{(k-1){\binom{m}{2}}}:({\rm SL}_m(2)\times \Sym_k)$.
\ep
\proof\ By Lemma \ref{LO} (3) and (4), $H=O_2(H):K$.
By Proposition \ref{PAutw}, $\Aut(R^k,w^k)=O(R,q)\wr\Sym_k$.

Clearly, $\Sym_k$ preserves $\mathcal{S}(\Phi,\Psi;k)$.
By Lemma \ref{LO} (3) and (\ref{EqT}), $O_2(H)_{(1i)}$ $(2\le i\le k)$ and $K_{(12\dots k)}$ also preserve $\mathcal{S}(\Phi,\Psi;k)$.
Hence $G$ contains $\langle O_2(H)_{(1i)}, K_{(12\dots k)},\Sym_k\mid 2\le i\le k\rangle$.

Let $g\in G$.
Up to $\Sym_k$, we may assume $g=(g_1,g_2,\dots,g_k)\in O(R,q)^{k}\cap G$.
If follows from $k\ge3$ that $g$ preserves $\Phi_{(1i)}$ for all $i$.
Hence $g_i\in H$ and $g_i\in O_2(H)g_1$.
Up to $K_{(12\dots k)}$, we may assume $g_i\in O_2(H)$ for all $i$ since $H=O_2(H):K$.
Moreover, up to $\langle O_2(H)_{(1i)}\mid 2\le i\le k\rangle$, we may assume $g_i=1$ if $i\ge2$.
Let $b\in\Psi$.
Then $g_1(b)=b+c$ for some $c\in\Phi$, and $(b,b,\dots,b)+g(b,b,\dots,b)=(c,0,\dots,0)\in \mathcal{S}(\Phi,\Psi;k)$.
Since $\Phi\cap\Psi=0$, we have $c=0$ by (\ref{EqT}), and $g_1=1$ by Lemma \ref{LO} (3).
Thus $g\in\langle O_2(H)_{(1i)}, K_{(12\dots k)},\Sym_k\mid 2\le i\le k\rangle$.

One can easily see that $G$ has the desired shape by Lemma \ref{LO} (3) and (4).
\qe

Later, we need the following lemma.

\bl\label{LN4} Let $v$ be a vector in $\mathcal{S}(\Phi,\Psi;3)$ with $w^3(v)=4$.
Assume that $\Phi\cap\Psi=0$.
Then one of the following holds:
\begin{enumerate}[{\rm (I)}]
\item $v=\sigma(a,a,0)$ for some $a\in\Phi\setminus\{0\}$ and $\sigma\in\Sym_3$.
\item $v=\sigma(a+c,a+b+c,b+c)$ for some $\sigma\in\Sym_3$, $a,b\in\Phi$ and $c\in\Psi\setminus\{0\}$ satisfying $w(a+c)=2$ and $w(a+b+c)=w(b+c)=1$. 
\end{enumerate}
Moreover, the numbers of vectors of type (I) and (II) are $3\times (2^m-1)$ and $3\times (2^m-1)\times 2^{2m-2}$ respectively.
\el
\proof\ The first assertion follows from (\ref{EqT}).
The number of vectors of type (I) can be calculated by $|\Phi\setminus\{0\}|=2^{m}-1$.
Let us count the number of vectors of type (II).
Note that $w(a+c)=2$ and $w(a+b+c)=w(b+c)=1$ if and only if $\langle a,c\rangle=0$ and $\langle b,c\rangle=1$.
Clearly, there are $2^{m}-1$ choices of $c\in\Psi\setminus\{0\}$.
For each $c$, there are $2^{m-1}$ choices of $a$ since $\{u\in \Phi\mid\langle u,c\rangle=0\}$ is an $(m-1)$-dimensional subspace of $\Phi$.
Similarly, there are $2^{m-1}$ choices of $b$ since $|\{u\in \Phi\mid\langle u,c\rangle=1\}|=2^m-2^{m-1}=2^{m-1}$.
Hence the number of vectors of type (II) is $3\times (2^{m}-1)\times 2^{2m-2}$.\qe

\subsection{Uniqueness of maximal totally singular subspaces for $k=3$}
In this subsection, we consider the case $k=3$, and show that any maximal totally singular subspace of $R^3$ satisfying (\ref{Cond1}) is conjugate to $\mathcal{S}(\Phi,\Psi;3)$ with $\Phi\cap\Psi=0$ under $O(R,q)^{3}$.
For the definition of $\mathcal{S}(\Phi,\Psi;3)$, see (\ref{EqT}).

Let $\rho_i$ denote the $i$-th coordinate projection from $R^3$ to $R$.
For a subspace $\mathcal{S}$ of $R^3$ and distinct $i,j\in \{1,2,3\}$, we denote $\mathcal{S}^{(i)}=\{v\in \mathcal{S}\mid \rho_i(v)=0\}$ and $\mathcal{S}^{(ij)}=\{v\in \mathcal{S}\mid \rho_i(v)=\rho_j(v)=0\}$.

\bl\label{LLem} Let $\mathcal{S}$ be a maximal totally singular subspace of $R^3$ satisfying (\ref{Cond1}) and let $i,j\in\{1,2,3\}$ such that $i\neq j$.
Then the following hold:
\begin{enumerate}
\item $\mathcal{S}^{(ij)}=0$.
\item $\rho_i(\mathcal{S})=R$.
\item The dimension of $\mathcal{S}^{(i)}$ is $m$.
\item $\rho_j(\mathcal{S}^{(i)})$ is a maximal totally singular subspace of $R$.
\end{enumerate}
\el
\proof\ Without loss of generality, we may assume that $i=1$ and $j=2$.

Let $v=(0,0,a)\in\mathcal{S}^{(12)}$.
Then $v$ must be $0$ by (\ref{Cond1}), which shows (1).

Suppose $\rho_1(\mathcal{S})\neq R$.
Then $\rho_1(\mathcal{S})^\perp\neq0$.
Let $a\in\rho_1(\mathcal{S})^\perp\setminus\{0\}$.
Then $(a,0,0)\in \mathcal{S}^\perp=\mathcal{S}$ by Lemma \ref{LO} (1), which contradicts (1).
Hence (2) holds.

By (2), the projection $\rho_1:\mathcal{S}\to R$ is surjective.
It follows from $\dim R=2m$ and $\dim \mathcal{S}=3m$ that $\dim \mathcal{S}^{(1)}=m$.
Hence we have (3).

Let $v=(0,a,b)\in \mathcal{S}^{(1)}$.
Since $v$ is singular, both $a$ and $b$ must be singular by (\ref{Cond1}).
Hence $\rho_2(\mathcal{S}^{(1)})$ is totally singular.
By (1) and (3), $\dim \rho_2(\mathcal{S}^{(1)})=m$, and hence $\rho_2(\mathcal{S}^{(1)})$ is a maximal totally singular subspace of $R$, which proves (4).\qe

Recall that $\Aut(R^3,w^3)\cong O(R,q)\wr\Sym_3$ from Proposition \ref{PAutw} and that $O(R,q)^3$ is the normal subgroup of $\Aut(R^3,w^3)$ isomorphic to the direct product of three copies of $O(R,q)$.

\bt\label{TCh} Let $\mathcal{S}$ be a maximal totally singular subspace of $(R^3,q^3)$ satisfying (\ref{Cond1}).
Then $\mathcal{S}$ is conjugate to $\mathcal{S}(\Phi,\Psi;3)$ under $O(R,q)^{3}$ for some maximal totally singular subspaces $\Phi$ and $\Psi$ of $(R,q)$ satisfying $\Phi\cap\Psi=0$.
Moreover $O(R,q)^{3}$ is transitive on the set of all maximal totally singular subspaces of $(R^3,q^3)$ satisfying (\ref{Cond1}).
\et
\proof\ 
Set $\Phi=\rho_1(\mathcal{S}^{(3)})$ and $\Phi^\prime=\rho_2(\mathcal{S}^{(3)})$.
Then both $\Phi$ and $\Phi'$ are maximal totally singular subspaces of $R$ by Lemma \ref{LLem} (4).
Up to the action of $O(R,q)$ on the second coordinate, we may assume that $\Phi'=\Phi$ by Lemma \ref{LO} (2).
It follows from Lemma \ref{LO} (3) that ${\rm Stab}_{O(R,q)}(\Phi')/O_2({\rm Stab}_{O(R,q)}(\Phi'))\cong {\rm SL}_m(2)$ acts naturally on $\Phi'$.
Hence we may assume that $\mathcal{S}^{(3)}=\{(a,a,0)\mid a\in \Phi\}$ up to the action of $O(R,q)$ on the second coordinate.
By the same arguments on $\rho_1(\mathcal{S}^{(2)})$ and $\rho_3(\mathcal{S}^{(2)})$, we obtain $\mathcal{S}^{(2)}=\{(a,0,a)\mid a\in \Phi''\}$ for some maximal totally singular subspace $\Phi''$ of $R$ up to the action of $O(R,q)$ on the third coordinate.
Since $\mathcal{S}$ is totally singular, we have $\langle\mathcal{S}^{(3)},\mathcal{S}^{(2)}\rangle=0$, equivalently, $\Phi''\subset\Phi^\perp$.
By Lemmas \ref{LO} (1) and \ref{LLem} (4), $\Phi''=\Phi$.
Hence we obtain two $m$-dimensional totally singular subspaces $\Phi_{(12)}=\{(a,a,0)\mid a\in \Phi\}$ and $\Phi_{(13)}=\{(a,0,a)\mid a\in \Phi\}$ of $\mathcal{S}$.

Let $c\in R\setminus \Phi$.
Then by Lemma \ref{LLem} (2), there exist $a,b\in R$ such that $v=(a,b,c)\in\mathcal{S}$.
It follows from $\langle\Phi_{(13)},v\rangle=0$ that $\langle p,a+\Phi\rangle=\langle p,c+\Phi\rangle$ for all $p\in \Phi$.
Note that $\{\langle p,\cdot\rangle\mid p\in \Phi\}=\Hom(R/\Phi,\F_2)$ by Lemma \ref{LO} (1).
Hence $a\in c+\Phi$.
Similarly, $a\in b+\Phi$.
Hence there exist $u\in \Phi_{(12)}$ and $x(c)\in \Phi$ such that $v+u=(x(c)+c,c,c)\in\mathcal{S}$. 
By (\ref{Cond1}) $x(c)$ is uniquely determined by $c$.
Hence for $c_1,c_2\in R\setminus \Phi$ with $c_1+c_2\notin \Phi$ we have $x(c_1)+x(c_2)=x(c_1+c_2)$.

Let $\Psi$ be a maximal totally singular subspace of $R$ such that $\Phi\cap \Psi=0$.
Set $\mathcal{T}=\{(x(c)+c,c,c)\mid c\in \Psi\}$.
Then $\mathcal{T}$ is an $m$-dimensional totally singular subspace of $\mathcal{S}$.
Since $\rho_2(\mathcal{T})=\rho_3(\mathcal{T})={\Psi}$ is totally singular, so is $\rho_1(\mathcal{T})$.
It follows from $\Phi\cap \Psi=0$ that $\Phi\cap \rho_1(\mathcal{T})=0$.
By Lemma \ref{LO} (5), $\rho_1(\mathcal{T})$ is conjugate to $\Psi$ under $O(R,q)$.
Hence we may assume $x(c)=0$ for all $c\in \Psi$ up to the action of $O(R,q)$ on the first coordinate.
Then $\mathcal{T}=\Psi_{(123)}=\{(c,c,c)\mid c\in\Psi\}$.
Hence $\mathcal{S}$ is conjugate to $\mathcal{S}(\Phi,\Psi;3)$ under $O(R,q)^{3}$.

The transitivity of $O(R,q)^3$ follows from Lemma \ref{LO} (5).\qe

\section{Application to binary codes and lattices}
In this section, we apply the results in the previous section to binary codes and lattices.
For the definitions and fundamental facts on binary codes and lattices, we refer the reader to \cite{MS,CS}.

\subsection{Application to binary codes}
In this subsection, we construct doubly even self-dual binary codes without codewords of weight $4$, which is a slight generalization of Turyn's construction of the extended binary Golay code (cf.\ \cite[Ch 18, Section 7.4]{MS}).

Let $(\ ,\ )$ be the standard inner product on $\F_2^n$.
For $x=(x_1,x_2,\dots,x_n)\in \F_2^n$, $\wt(x)=|\{i\mid x_i\neq0\}|$ is the (Hamming) {\it weight} of $x$.
For $x,y\in\F_2^n$, the following holds:
\begin{eqnarray}
( x,y)\equiv\frac{1}{2}(\wt(x+y)-\wt(x)-\wt(y))\pmod 2\label{Eq:qC}.
\end{eqnarray}
A subset of $\F_2^n$ is called a {\it binary ({\it linear}) code} of {\it length} $n$ if it is a subspace.
Let $C^\perp$ denote the dual code of a binary code $C$ of length $n$, that is, $C^\perp=\{c\in \F_2^n\mid ( c, C)=0\}$.
A binary code $C$ is said to be {\it doubly even} if $\wt(c)\in4\Z$ for all $c\in C$, and is said to be {\it self-dual} if $C=C^\perp$.
We now consider the following condition on $C^\perp$:
\begin{eqnarray}
\wt(c)\in2\Z\quad \forall c\in C^\perp.\label{CondC}
\end{eqnarray}
Note that $C^\perp$ satisfies (\ref{CondC}) if and only if $C$ contains the all-one codeword $\1_n$.

Let $C$ be a doubly even binary code of length $n\in8\Z$ satisfying (\ref{CondC}).
Set $R(C)=C^\perp/C$ and let $$\varphi_C:C^\perp\to R(C),\quad c\mapsto c+C$$ be the canonical map.
Then $R(C)\cong (\Z/2\Z)^m$ for some $m\in\Z_{\ge0}$ and we view it as a vector space over $\F_2$.
Note that $\dim C=(n-m)/2$.
Let $\langle\ ,\ \rangle$ be the symmetric bilinear form on $R(C)$ defined by $\langle c+C,d+C\rangle=(c,d)$.
By (\ref{CondC}) and $(C^\perp)^\perp=C$, it is a non-singular symplectic form.
Consider the map $$q_C:R(C)\to\F_2,\quad x+C\mapsto \wt (x)/2 \pmod2.$$
It follows from $\wt(c)\in4\Z$ for all $c\in C$, $(C,C^\perp)=0$ and (\ref{Eq:qC}) that $q_C$ is a well-defined quadratic form on $R(C)$ associated to the symplectic form $\langle\ ,\ \rangle$.
Since $n\in8\Z$, the type of $q_C$ is plus.
Thus $(R(C),q_C)$ is a non-singular $m$-dimensional quadratic space of plus type over $\F_2$.

Let $\mathcal{S}$ be a maximal totally singular subspace of $(R(C)^k,q_C^k)$.
Let $\mathfrak{C}(\mathcal{S})$ be the inverse image of ${\mathcal{S}}$ with respect to the canonical map 
\begin{eqnarray*}
\varphi^k_C:(C^\perp)^k\to R(C)^k,\quad (c_i)\mapsto(\varphi_C(c_i)),\label{Eq:mapk}
\end{eqnarray*}
that is, $$\mathfrak{C}(\mathcal{S})=(\varphi^{k}_C)^{-1}(\mathcal{S}).$$
Then $\mathfrak{C}(\mathcal{S})$ is a binary code of length $nk$.
Note that $\Ker\varphi^k_C=C^k$, where $C^k$ is the direct sum of $k$ copies of $C$, and that $\dim \Ker\varphi^k_C=(n-m)k/2$.

\bt\label{MTC} Let $C$ be a doubly even binary code of length $n\in8\Z$ satisfying (\ref{CondC}) and let $\mathcal{S}$ be a maximal totally singular subspace of $(R(C)^k,q_C^k)$.
\begin{enumerate}[{\rm (1)}]
\item The binary code $\mathfrak{C}(\mathcal{S})$ of length $nk$ is doubly even and self-dual.
\item Assume that $\mathcal{S}$ satisfies (\ref{Cond1}).
If $C$ has no codewords of weight $4$ then so does $\mathfrak{C}(\mathcal{S})$.
\end{enumerate}
\et
\proof\ Since ${\mathcal{S}}$ is totally singular, $\mathfrak{C}(\mathcal{S})$ is doubly even.
It follows from the maximality of $\mathcal{S}$ that $\dim \mathcal{S}=mk/2$.
Hence $$\dim\mathfrak{C}(\mathcal{S})=\dim\Ker\varphi_C^k+\dim \mathcal{S}=k(n-m)/2+mk/2=kn/2,$$ which shows that $\mathfrak{C}(\mathcal{S})$ is self-dual.

By the definition of the map $w: R(C)\to\{0,1,2\}$ (cf.\ Section 2.1), the minimum weight of $c+C\in R(C)$ is greater than or equal to $2w(c+C)$.
Hence for $u\in R(C)^k$ the minimum weight of $(\varphi_C^k)^{-1}(u)$ is greater than or equal to $2w^k(u)$.
Thus (2) follows from (\ref{Cond1}) and the assumption on $C$.
\qe

\br In general, the minimum weight of $c+C\in R(C)$ may not be equal to $2w(c+C)$.
\er

\br Let $C$ be a doubly even binary code of length $n\in8\Z$ satisfying (\ref{CondC}).
Let $A$ and $B$ be doubly even self-dual binary codes of length $n$ satisfying $A\cap B\supset C$.
Then $\Phi=\varphi_C(A)$ and $\Psi=\varphi_C(B)$ are maximal totally singular subspaces of $R(C)$.
By (\ref{EqT}) $\mathfrak{C}(\mathcal{S}(\Phi,\Psi;k))$ is described in terms of $A$ and $B$ as follows:
$$\mathfrak{C}(\mathcal{S}(\Phi,\Psi;k))= \{(a_1+b,a_2+b,\dots,a_{k}+b)\in\F_2^{nk}\mid a_i\in A,b\in B, \sum_{i=1}^k a_i\in  A\cap B\}.$$

\er

\subsection{Turyn's construction of the extended binary Golay code}
In this subsection, we consider the case where $k=3$ and $C={\rm Span}_{\F_2}\{\1_8\}$, and describe the extended binary Golay code of length $24$.
Moreover, as an application of Section 2.2, we study the extended binary Golay code and the Mathieu group.
We continue the notation of the previous subsection.

Obviously, $C$ has no codewords of weight $4$ and satisfies (\ref{CondC}).
Let $\mathcal{S}$ be a maximal totally singular subspace of $(R(C)^3,q_C^3)$ satisfying (\ref{Cond1}).
Note that such a subspace exists by Proposition \ref{PLM}.
By Theorem \ref{MTC}, $\mathfrak{C}(\mathcal{S})$ is a doubly even self-dual binary code of length $24$ without codewords of weight $4$. 
It is well-known that such a binary code is equivalent to the extended binary Golay code of length $24$ (cf.\ \cite{MS}).

\bc Let $C={\rm Span}_{\F_2}\{\1_8\}$ and let $\mathcal{S}$ be a maximal totally singular subspace of $(R(C)^3,q_C^3)$ satisfying (\ref{Cond1}).
Then $\mathfrak{C}(\mathcal{S})$ is equivalent to the extended binary Golay code of length $24$.
\ec

\br Let $\Phi$ and $\Psi$ be maximal totally singular subspaces of $(R(C),q_C)$ satisfying $\Phi\cap\Psi=0$.
Then $\mathfrak{C}(\mathcal{S}(\Phi,\Psi;3))$ is isomorphic to the extended binary Golay code, which is called Turyn's construction {\rm (cf.\ \cite[Ch 18, Theorem 12]{MS})}.
Note that $\varphi_C^{-1}(\Phi)$ and $\varphi_C^{-1}(\Psi)$ are isomorphic to the extended Hamming code of length $8$.
\er

Since $\Aut(C^\perp)=\Aut(C)$, there is a canonical homomorphism of groups
\begin{eqnarray*}
\psi_C: \Aut(C^\perp)\to O(R(C),q_C).\label{Eq:psiC}
\end{eqnarray*}
The injectivity is clear, and the surjectivity follows from $\Aut(C)=\Sym_8\cong O(R(C),q_C)$.
Hence $\psi_C$ is an isomorphism of groups. 

A subcode of $\mathfrak{C}(\mathcal{S})$ is called a {\it trio} if it is equivalent to $C^3=C\oplus C\oplus C$.
Since $\Aut(C^3)\cong \Aut(C)\wr\Sym_3$ and $\Aut(C^3)=\Aut((C^\perp)^3)$, $\psi_C$ induces the injective homomorphism of groups
\begin{eqnarray*}
\psi^3_C:\Aut((C^\perp)^3)\to\Aut(R(C)^3,w^3).\label{EqCCC}
\end{eqnarray*}
By Proposition \ref{PAutw}, $\psi^3_C$ is an isomorphism.
Note that $\psi^3_C$ is compatible with $\varphi_C^3$, that is, $\varphi^3_C(g(c))=\psi^3_C(g)(\varphi^3_C(c))$ for $g\in\Aut((C^\perp)^3)$ and $c\in (C^\perp)^3$.

As corollaries of Proposition \ref{PStab} and Theorem \ref{TCh}, we obtain the following well-known properties of $\Aut(\mathfrak{C}(\mathcal{S}))$. 

\bc\label{CMTC} {\rm (cf.\ \cite{CS})} Let $C={\rm Span}_{\F_2}\{\1_8\}$ and let $\mathcal{S}$ be a maximal totally singular subspace of $(R(C)^3,q_C^3)$ satisfying (\ref{Cond1}).
Let $G$ be the stabilizer of $C^3$ in the automorphism group of $\mathfrak{C}(\mathcal{S})$.
\begin{enumerate}[{\rm (1)}]
\item The automorphism group of $\mathfrak{C}(\mathcal{S})$ is transitive on the set of all trios of $\mathfrak{C}(\mathcal{S})$.
\item $G$ is isomorphic to the stabilizer of $\mathcal{S}$ in $\Aut(R(C)^3,w^3)$, and it has the shape $2^{6}:({\rm SL}_3(2)\times\Sym_3)$.

\end{enumerate}
\ec
\proof\ Let $D$ be a trio of $\mathfrak{C}(\mathcal{S})$.
It suffices to show that there is an automorphism of $\mathfrak{C}(\mathcal{S})$ which sends $D$ to $C^3$.
Let $g\in\Sym_{24}$ such that $g(D)=C^3$.
Since $g(\mathfrak{C}(\mathcal{S}))$ is a doubly even self-dual binary code without codewords of weight $4$,  its image $\varphi^3_C(g(\mathfrak{C}(\mathcal{S})))$ is a maximal totally singular subspace of $R(C)^3$ satisfying (\ref{Cond1}).
Then by Theorem \ref{TCh}, there exists $h\in \Aut(C^3)(\cong \Aut(R(C)^3,w^3))$ such that $h(C^3)=C^3$ and $h( g(\mathfrak{C}(\mathcal{S})))=\mathfrak{C}(\mathcal{S})$, which proves (1).

Clearly, $G$ is equal to the stabilizer of $\mathfrak{C}(\mathcal{S})$ in $\Aut((C^\perp)^3)$.
Since $\psi_C^3$ is compatible with $\varphi^3_C$ and $\varphi^3_C(\mathfrak{C}(\mathcal{S}))=\mathcal{S}$, $\psi_C^3(G)$ is the stabilizer of $\mathcal{S}$ in $\Aut(R(C)^3,w^3)$.
Since $\psi_C^3$ is an isomorphism, $G$ is isomorphic to the stabilizer of $\mathcal{S}$ in $\Aut(R(C)^3,w^3)$.
Moreover, $G$ is isomorphic to the stabilizer of $\mathcal{S}(\Phi,\Psi;3)$ in $\Aut(R(C)^3,w^3)$ by (1).
Hence we obtain (2) by Proposition \ref{PStab}.\qe

Let us count the number of codewords of weight $8$ in $\mathfrak{C}(\mathcal{S})$.
The following lemma is easy.

\bl\label{LCc} Let $u\in R(C)$.
\begin{enumerate}
\item If $w(u)=0$ then the number of codewords of weight $8$ in $\varphi_C^{-1}(u)$ is $1$.
\item If $w(u)=1$ then the number of codewords of weight $2$ in $\varphi_C^{-1}(u)$ is $1$.
\item If $w(u)=2$ then the number of codewords of weight $4$ in $\varphi_C^{-1}(u)$ is $2$.
\end{enumerate}
\el

\bp Let $C={\rm Span}_{\F_2}\{\1_8\}$ and let $\mathcal{S}$ be a maximal totally singular subspace of $(R(C)^3,q_C^3)$ satisfying (\ref{Cond1}).
Then the number of codewords of weight $8$ in $\mathfrak{C}(\mathcal{S})$ is $759$.
\ep
\proof\ By Theorem \ref{TCh} and the surjectivity of $\psi_C^3$, we may assume that $\mathcal{S}=\mathcal{S}(\Phi,\Psi;3)$.
Any codeword of weight $8$ in $\mathfrak{C}(\mathcal{S})$ belongs to $(\varphi_C^{3})^{-1}(0)$ or $(\varphi_C^{3})^{-1}(u)$ for some $u\in R(C)^3$ with $w^3(u)=4$.
Note that $(\varphi_C^{3})^{-1}(0)=C^3$.
Hence by Lemmas \ref{LN4} and \ref{LCc}, the number of codewords of weight $8$ in $\mathfrak{C}(\mathcal{S})$ is $$ (3\times 1)+(3\times (2^3-1))\times 1^0\times2^2+(3\times (2^3-1)\times 2^4)\times 1^2\times 2^1=759.$$\qe

\subsection{Application to lattices}
In this subsection, we construct even unimodular lattices without vectors of norm $2$, which was studied by Griess \cite{Gr}.

Let $(\ , \ )$ be a positive definite symmetric bilinear form on $\R^n$.
For $x\in \R^n$, $( x,x)$ is the (squared) {\it norm} of $x$.
For $x,y\in\R^n$, the following holds:
\begin{eqnarray}
2( x, y)= ( x+y,x+y)-( x,x)-( y,y)\label{Eq:qL}.
\end{eqnarray}
A subset $L$ of $\R^n$ is called a {\it lattice} of {\it rank} $n$ if $L$ has a basis $e_1,e_2,\dots,e_n$ of $\R^n$ satisfying $L=\oplus_{i=1}^n\Z e_i$.
Let $L^*$ denote the dual lattice of a lattice $L$ of rank $n$, that is, $L^*=\{v\in \R^n\mid ( v, L)\subset\Z\}$.
A lattice $L$ is said to be {\it even} if $( v,v)\in2\Z$ for all $v\in L$, and is said to be {\it unimodular} if $L=L^*$.
We now consider the following condition on $L^*$:
\begin{eqnarray}
( v,v)\in\Z\quad \forall v\in L^*.\label{CondL}
\end{eqnarray}

\bl\label{LL1} Let $L$ be an even lattice of rank $n\in8\Z$.
\begin{enumerate}[{\rm (1)}]
\item $L$ satisfies (\ref{CondL}) if and only if $L$ contains $\sqrt2 J$ for some even unimodular lattice $J$.
\item If $L$ satisfies (\ref{CondL}) then $2L^*\subset L$.
\end{enumerate}
\el
\proof\ If $L\supset \sqrt2 J$ for some even unimodular lattice $J$ then $L^*\subset J/\sqrt2$.
Since $( v,v)\in\Z$ for any $v\in J/\sqrt2$, $L^*$ satisfies (\ref{CondL}).

Conversely, we assume that $L$ satisfies (\ref{CondL}).
Then $\sqrt2L^*$ is even.
Since $n\in8\Z$, there exists an even unimodular lattice $J$ of rank $n$ such that $\sqrt2L^*\subset J$ (cf.\ \cite{Ve}).
Hence $J\subset (\sqrt2L^*)^*=L/\sqrt2$, equivalently, $\sqrt2J\subset L$, and obtain (1).

If $L$ satisfies (\ref{CondL}) then $\sqrt2 L^*\subset(\sqrt2 L^*)^*= L/\sqrt2$, and obtain (2).
\qe

Let $L$ be an even lattice of rank $n\in8\Z$ satisfying (\ref{CondL}).
Set $R(L)=L^*/L$ and let $$\varphi_L:L^*\to R(L),\quad v\mapsto v+L$$ be the canonical map.
Then by Lemma \ref{LL1} (2), $R(L)\cong (\Z/2\Z)^m$ for some $m\in\Z_{\ge0}$ and we view it as a vector space over $\F_2$.
Note that the determinant of $L$ is $2^{m/2}$.
Let $\langle\ ,\ \rangle$ be the symmetric bilinear form on $R(L)$ defined by $\langle v+L,u+L\rangle=2(v,u)\pmod2$.
By (\ref{CondL}) and $(L^*)^*=L$, it is a non-singular symplectic form.
Consider the map $$q_L:R(L)\to\F_2,\quad v+L\mapsto ( v,v)\pmod 2.$$
It follows from $(v,v)\in2\Z$ for all $v\in L$, $( L,L^*)\subset \Z$ and (\ref{Eq:qL}) that $q_L$ is a well-defined quadratic form on $R(L)$ associated to $\langle\ ,\ \rangle$.
Since $n\in8\Z$, the type of $q_L$ is plus (cf.\ \cite{Ve}).
Thus $(R(L),q_L)$ is a non-singular $m$-dimensional quadratic space of plus type over $\F_2$.

Let $\mathcal{S}$ be a maximal totally singular subspace of $(R(L)^k,q_L^k)$.
Let $\mathfrak{L}(\mathcal{S})$ be the inverse image of ${\mathcal{S}}$ with respect to the canonical map 
\begin{eqnarray*}
\varphi^k_L:(L^*)^k\to R(L)^k,\quad (v_i)\mapsto(\varphi_L(v_i)),\label{Eq:mapkL}
\end{eqnarray*}
that is $$\mathfrak{L}(\mathcal{S})=(\varphi_L^k)^{-1}({\mathcal{S}}).$$
Then $\mathfrak{L}(\mathcal{S})$ is a lattice, and its rank is $nk$ since it contains $L^k$, where $L^k$ is the orthogonal direct sum of $k$ copies of $L$.

\bt\label{MTL} {\rm (cf.\ \cite{Gr})} Let $L$ be an even lattice satisfying (\ref{CondL}) and let $\mathcal{S}$ be a maximal totally singular subspace of $(R(L)^k,q_L^k)$.
\begin{enumerate}
\item
The lattice $\mathfrak{L}(\mathcal{S})$ of rank $nk$ is even and unimodular.
\item Assume that $\mathcal{S}$ satisfies (\ref{Cond1}).
If $L$ has no vectors of norm $2$ then so does $\mathfrak{L}(\mathcal{S})$.
\end{enumerate}
\et
\proof\ Since ${\mathcal{S}}$ is totally singular, $\mathfrak{L}(\mathcal{S})$ is even.
By the maximality of $\mathcal{S}$, we have $\dim \mathcal{S}=mk/2$.
Hence the determinant of ${\mathfrak{L}}_L(\mathcal{S})$ is $2^{mk/2-mk/2}=1$, which shows that $\mathfrak{L}(\mathcal{S})$ is unimodular.

By the definition of the map $w:R(L)\to\{0,1,2\}$, the minimum norm of $v+L\in R(L)$ is greater than or equal to $w(v+L)$.
Hence for $u\in R(L)^k$ the minimum norm of $(\varphi_L^k)^{-1}(u)$ is greater than or equal to $w^k(u)$.
Thus (2) follows from (\ref{Cond1}) and the assumption on $L$.
\qe

\br In general, the minimum norm of $\lambda+L\in R(L)$ may not be equal to $w(\lambda+L)$.
\er

\br Let $L$ be an even lattice of rank $n\in8\Z$ satisfying (\ref{CondL}).
Let $J$ and $K$ be even unimodular lattices of rank $n$ satisfying $J\cap K\supset L$.
Then $\Phi=\varphi_L(J)$ and $\Psi=\varphi_L(K)$ are maximal totally singular subspaces of $R(L)$.
By (\ref{EqT}) $\mathfrak{L}(\mathcal{S}(\Phi,\Psi;k))$ is described in terms of $J$ and $K$ as follows:
$$\mathfrak{L}(\mathcal{S}(\Phi,\Psi;k))=\{(a_1+b,a_2+b,\dots,a_{k}+b)\in\R^{nk}\mid a_i\in J,b\in K,\sum_{i=1}^ka_i\in J\cap K\}.$$
\er

\bn In \cite{Gr}, $\mathfrak{L}(\mathcal{S}(\Phi,\Psi;k))$ was constructed and its minimum norm was studied.
\en

\subsection{Lepowsky and Meurman's description of the Leech lattice}
Let ${E_8}$ denote the $E_8$ root lattice.
In this subsection, we consider the case where $k=3$ and $L=\sqrt2{E_8}$, and describe the Leech lattice.
Moreover, as an application of Section 2.2, we study the Leech lattice and the Conway group.
Note that the method in this subsection was already considered in \cite{LM,Gr2}.
We continue the notation of the previous subsection.

Since $E_8$ is even, $L$ has no vectors of norm $2$.
By Lemma \ref{LL1} (1), $L$ satisfies (\ref{CondL}).
Let $\mathcal{S}$ be a maximal totally singular subspace of $(R(L)^3,q_L^3)$ satisfying (\ref{Cond1}).
Note that such a subspace exists by Proposition \ref{PLM}.
By Theorem \ref{MTL}, $\mathfrak{L}(\mathcal{S})$ is an even unimodular lattice of rank $24$ without vectors of norm $2$. 
It was shown in \cite{Co} that such a lattice is isomorphic to the Leech lattice.

\bc\label{CL} Let $L=\sqrt2E_8$ and let $\mathcal{S}$ be a maximal totally singular subspace of $(R(L)^3,q_L^3)$ satisfying (\ref{Cond1}).
Then $\mathfrak{L}(\mathcal{S})$ is isomorphic to the Leech lattice.
\ec

\br Let $\Phi$ and $\Psi$ be maximal totally singular subspaces of $(R(L),q_L)$ satisfying $\Phi\cap\Psi=0$.
Then $\mathfrak{L}(\mathcal{S}(\Phi,\Psi;3))$ is isomorphic to the Leech lattice, which was shown in \cite[Corollary 2.3]{LM}).
Note that $\varphi_L^{-1}(\Phi)$ and $\varphi_L^{-1}(\Psi)$ are isomorphic to the $E_8$ root lattice.
\er

Since $\Aut(L^*)=\Aut(L)$, there is a canonical homomorphism of groups
\begin{eqnarray*}
\psi_L: \Aut(L^*)\to O(R(L),q_L).\label{Eq:psiL}
\end{eqnarray*}
Moreover, the center $Z(\Aut(L))$ of $\Aut(L)$ is generated by the $-1$-isometry, and the quotient group $\Aut(L)/Z(\Aut(L))$ is isomorphic to $O(R(L),q_L)$ (cf.\ \cite{CS}).
Hence $\psi_L$ is surjective and $\Ker\psi_L=Z(\Aut(L))$. 

Since the automorphism group of $L^3$ is isomorphic to $\Aut(L)\wr\Sym_3$ and $\Aut(L^3)=\Aut((L^*)^3)$, 
$\psi_L$ induces a homomorphism of groups
\begin{eqnarray*}
\psi^3_L:\Aut((L^*)^3)\to\Aut(R(L)^3,w^3).\label{EqLLL}
\end{eqnarray*}
By Proposition \ref{PAutw}, $\psi^3_L$ is surjective.
It follows from $\Ker\psi_L=Z(\Aut(L))$ that $\Ker\psi^3_L=Z(\Aut(L))^{3}\cong 2^3$.
Note that $\psi^3_L$ is compatible with $\varphi_L^3$, that is, $\varphi^3_L(g(v))=\psi^3_L(g)(\varphi^3_L(v))$ for $g\in\Aut((L^*)^3)$ and $v\in (L^*)^3$.

As corollaries of Proposition \ref{PStab} and Theorem \ref{TCh}, we obtain the following well-known properties of $\Aut(\mathfrak{L}(\mathcal{S}))$. 

\bc\label{MCL} {\rm \cite{Gr2}} Let $L=\sqrt2E_8$ and let $\mathcal{S}$ be a maximal totally singular subspace of $(R(L)^3,q_L^3)$ satisfying (\ref{Cond1}).
Let $G$ be the stabilizer of $L^3$ in the automorphism group of $\mathfrak{L}(\mathcal{S})$.
\begin{enumerate}[{\rm (1)}]
\item The automorphism group of $\mathfrak{L}(\mathcal{S})$ is transitive on the set of all sublattices of $\mathfrak{L}(\mathcal{S})$ isomorphic to the orthogonal direct sum of three copies of $\sqrt2E_8$.
\item $G$ is isomorphic to an extension of the stabilizer of $\mathcal{S}$ in $\Aut(R(L)^3,w^3)$ by an elementary abelian $2$-group of order $2^3$, and $G$ has the shape $2^3.(2^{12}:({\rm SL}_4(2)\times\Sym_3))$.
\end{enumerate}
\ec
\proof\ Let $J$ be a sublattice of $\mathfrak{L}(\mathcal{S})$ isomorphic to $L^{3}$.
It suffices to show that there is an automorphism of $\mathfrak{L}(\mathcal{S})$ which sends $J$ to $L^3$.
Let $g$ be an orthogonal transformation of $\R^{24}$ such that $g(J)=L^3$.
Then $g(\mathfrak{L}(\mathcal{S}))$ is an even unimodular lattice of rank $24$ without vectors of norm $2$ and contains $L^3$.
Hence $\varphi^3_L(g(\mathfrak{L}(\mathcal{S})))$ is a maximal totally singular subspace of $R(L)^3$ satisfying (\ref{Cond1}).
Then by Theorem \ref{TCh}, there exists $h\in \Aut(R(L)^3,w^3)$ such that $h( \psi^3_L(g(\mathfrak{L}(\mathcal{S}))))=\mathcal{S}$.
Since $\varphi_L^3$ is surjective and compatible with $\psi_L^3$, there exists $\tilde{h}\in\Aut(L^3)$ such that $\tilde{h}\circ g(\mathfrak{L}(\mathcal{S}))=\mathfrak{L}(\mathcal{S})$.
Since $\tilde{h}\circ g(J)=L^3$, we have proved (1).

Clearly, $G$ is equal to the stabilizer of $\mathfrak{L}(\mathcal{S})$ in $\Aut((L^*)^3)$.
Since $\psi_L^3$ is compatible with $\varphi^3_L$ and $\varphi^3_L(\mathfrak{L}(\mathcal{S}))=\mathcal{S}$, $\psi_L^3(G)$ is the stabilizer of $\mathcal{S}$ in $\Aut(R(L)^3,w^3)$.
Moreover, $\psi_L^3(G)$ is isomorphic to the stabilizer of $\mathcal{S}(\Phi,\Psi;3)$ in $\Aut(R(L)^3,w^3)$ by (1).
Since $\Ker\psi_L^3\cong 2^3$, we obtain (2) by Proposition \ref{PStab}.\qe

\bn Corollary \ref{MCL} (1) and (2) were shown in \cite[Proposition 3.5]{Gr2} and \cite[Corollary 3.10]{Gr2} respectively by a similar approach.
\en

Let us count the number of vectors of norm $4$ in $\mathfrak{L}(\mathcal{S})$.
The following lemma is easy.

\bl\label{LCl} Let $u\in R(L)$.
\begin{enumerate}
\item If $w(u)=0$ then the number of vectors of norm $4$ in $\varphi_L^{-1}(u)$ is $240$.
\item If $w(u)=1$ then the number of vectors of norm $1$ in $\varphi_L^{-1}(u)$ is $2$.
\item If $w(u)=2$ then the number of vectors of norm $2$ in $\varphi_L^{-1}(u)$ is $16$.
\end{enumerate}
\el

\bp Let $L=\sqrt2E_8$ and let $\mathcal{S}$ be a maximal totally singular subspace of $(R(L)^3,q_L^3)$ satisfying (\ref{Cond1}).
Then the number of vectors of norm $4$ in $\mathfrak{L}(\mathcal{S})$ is $196560$.
\ep
\proof\ By Theorem \ref{TCh} and the surjectivity of $\psi_L^3$, we may assume that $\mathcal{S}=\mathcal{S}(\Phi,\Psi;3)$.
Any vector of norm $4$ in $\mathfrak{L}(\mathcal{S})$ belongs to $(\varphi_L^{3})^{-1}(0)$ or $(\varphi_L^{3})^{-1}(u)$ for some $u\in R(L)^3$ with $w^3(u)=4$.
Note that $(\varphi_L^{3})^{-1}(0)=L^3$.
Hence by Lemmas \ref{LN4} and \ref{LCl}, the number of vectors of norm $4$ in $\mathfrak{L}(\mathcal{S})$ is  $$(3\times240)+(3\times (2^4-1))\times 2^0\times16^2+(3\times (2^4-1)\times 2^6)\times 2^2\times 16^1=196560.$$\qe 

\section{Application to vertex operator algebras}
In this section, applying the results in Section 2 to VOAs, we obtain the moonshine VOA as a simple current extension of the tensor product of three copies of $V_{\sqrt2E_8}^+$ and describe some automorphism group of the moonshine VOA.

The main approach in this section is to deduce the argument of VOAs to that of quadratic spaces in Section 2, which is quite similar to the approaches in Section 3 for binary codes and lattices.
For example, the transitivity of the automorphism group of the moonshine VOA on the set of subVOAs isomorphic to $(V_{\sqrt2E_8}^+)^3$ will be shown in Theorem \ref{MT} by using the uniqueness of certain maximal totally singular subspaces in Theorem \ref{TCh}, which is similar to Corollaries \ref{CMTC} and \ref{MCL}.

For the definitions and facts on VOAs and modules, see Section 1.

\subsection{A construction of holomorphic VOAs}

Let $V$ be a simple rational $C_2$-cofinite VOA of CFT type of central charge $n$.
Let $R(V)$ denote the set of all isomorphism classes of irreducible $V$-modules.
We consider the following conditions:
\begin{enumerate}[(a)]
\item Any irreducible $V$-module is a self-dual simple current module.
\item For any $[M]\in R(V)\setminus\{[V]\}$, the lowest weight of $M$ belongs to $\Z_{>0}/2$.
\item Let $q_V:R(V)\to\F_2$ be the map defined by setting $q_V([M])=0$ if $M$ is $\Z$-graded, and $q_V([M])=1$ if $M$ is $(\Z+1/2)$-graded.
Then $(R(V),q_V)$ is a non-singular quadratic space of plus type over $\F_2$.
\item For any irreducible $V$-module, the invariant bilinear form on it is symmetric.
\end{enumerate}

\br The assumption (b) corresponds to (\ref{CondC}) in binary codes and (\ref{CondL}) in lattices.
\er

If $V$ satisfies (a) then by Lemma \ref{LGp} $R(V)$ is an elementary abelian $2$-group under the fusion rules.
In this case, we view $R(V)$ as an $m$-dimensional vector space over $\F_2$, where $|R(V)|=2^m$.

\bl\label{LSCE} {\rm \cite{DM0} (cf.\ \cite[Proposition 1.1]{Sh4})} Let $V$ be a simple rational $C_2$-cofinite VOA of CFT type satisfying (a).
Let $X$ be a simple VOA containing $V$ as a full subVOA.
Then the multiplicity of any irreducible $V$-submodule of $X$ is $1$, and $X$ is a simple current extension of $V$.
\el

\bl\label{LSCE2} {\rm \cite[Theorem 1]{LY2}} Let $V$ be a simple rational $C_2$-cofinite VOA of CFT type satisfying (a).
Let $X=\oplus_{[M]\in\mathcal{T}}M$ be a simple current extension of $V$, where $\mathcal{T}$ is a subgroup of $R(V)$.
Let $U$ be an irreducible $X$-module and let $N$ be an irreducible $V$-submodule of $U$.
Set $\mathcal{U}=\{W\times [N]\mid W\in\mathcal{T}\}$.
Then $U\cong \oplus_{[M]\in \mathcal{U}}M$ as $V$-modules.
Moreover, $\oplus_{[M]\in \mathcal{U}}M$ has a unique irreducible $X$-module structure extending its $V$-module structure up to isomorphism.
\el

If $V$ satisfies (a)--(c) then the form $\langle\ ,\ \rangle:R(V)\times R(V)\to\F_2$ defined by setting $$\langle W,W'\rangle=q_V(W)+q_V(W')+q_V(W\times W')$$ is non-singular and symplectic.

Let $V^k$ be the tensor product of $k$ copies of $V$.
Then $V^k$ is a simple rational $C_2$-cofinite VOA of CFT type of central charge $nk$ (\cite{ABD,DMZ,FHL}).
By Lemma \ref{LemFHL}, if $V$ satisfies (a)--(d) then so does $V^k$, and $(R(V)^k,q_V^k)=(R(V^k),q_{V^k})$.
Let $\mathcal{T}$ be a subset of $R(V)^k$ and let $\mathfrak{V}(\mathcal{T})$ denote the $V^{ k}$-module defined by  
\begin{eqnarray*}
\mathfrak{V}(\mathcal{T})=\bigoplus_{[M]\in{\mathcal{T}}}M,\label{Eq:mapkV}
\end{eqnarray*}
where we identify $([M(1)],[M(2)],\dots,[M(k)])\in R(V)^k$ with the isomorphism class of irreducible $V^{ k}$-module $\otimes_{i=1}^k M(i)$ by Lemma \ref{LemFHL}.

\bp\label{PMax} Let $V$ be a simple rational $C_2$-cofinite VOA of CFT type satisfying (a)--(d) and let $\mathcal{T}$ be a subset of $R(V)^k$.
Then $\mathfrak{V}(\mathcal{T})$ has a simple VOA structure extending its $V^k$-module structure if and only if $\mathcal{T}$ is a totally singular subspace.
Moreover, $\mathfrak{V}(\mathcal{T})$ is holomorphic if and only if $\mathcal{T}$ is a maximal totally singular subspace.
\ep
\proof\ Recall that $V^k$ is a simple rational $C_2$-cofinite VOA of CFT type satisfying (a)--(d).
By (c), $\mathfrak{V}(\mathcal{T})$ is $\Z$-graded if and only if $\mathcal{T}$ is totally singular.
Hence the first assertion follows from Proposition \ref{PVOA} and (d).

Assume that $\mathfrak{V}(\mathcal{T})$ is holomorphic.
Let $\mathcal{U}$ be a totally singular subspace of $R(V)^k$ containing $\mathcal{T}$.
Then by the first assertion and Proposition \ref{PSC}, the VOA $\mathfrak{V}(\mathcal{U})$ containing $\mathfrak{V}(\mathcal{T})$ as a full subVOA.
Note that both $\mathfrak{V}(\mathcal{U})$ and $\mathfrak{V}(\mathcal{T})$ are of CFT type by (b).
Since $\mathfrak{V}(\mathcal{T})$ is holomorphic, $\mathfrak{V}(\mathcal{U})=\mathfrak{V}(\mathcal{T})$, and hence $\mathcal{U}=\mathcal{T}$.

Assume that $\mathcal{T}$ is a maximal totally singular subspace.
Let $U$ be an irreducible $\mathfrak{V}(\mathcal{T})$-module.
Take an irreducible $V^k$-submodule $N$ of $U$.
By Lemma \ref{LSCE2}, $U\cong \oplus_{[M]\in \mathcal{U}}M$ as $V^k$-modules, where $\mathcal{U}=\{W\times [N]\mid W\in\mathcal{T}\}$.
Let $W'\in\mathcal{T}$.
Then $W'\times[N]\in\mathcal{U}$.
Since $U$ is irreducible, $q_V^k([N])=q_V^k( W'\times [N])$.
By $q_V^k(W')=0$ we have $\langle W',[N]\rangle=q_V^k(W')+q_V^k([N])+q_V^k(W'\times [N])=0$.
Since $\mathcal{T}$ is maximal, we have $\mathcal{T}^\perp=\mathcal{T}$ by Lemma \ref{LO} (1), and hence $[N]\in\mathcal{T}$.
Thus $\mathcal{T}=\mathcal{U}$, and $U\cong\mathfrak{V}(\mathcal{T})$ as $\mathfrak{V}(\mathcal{T})$-modules by Lemma \ref{LSCE2}.
\qe

\br By Proposition \ref{PSC}, the VOA structure on $\mathfrak{V}(\mathcal{T})$ extending its $V^k$-module structure is uniquely determined by $\mathcal{T}$ up to isomorphism.
\er

\bt\label{MTV} Let $V$ be a simple rational $C_2$-cofinite VOA of CFT type of central charge $n$ satisfying (a)--(d).
Let $\mathcal{S}$ be a maximal totally singular subspace of $(R(V)^k,q_V^k)$.
\begin{enumerate}
\item $V^{ k}$-module $\mathfrak{V}(\mathcal{S})$ has a holomorphic VOA structure of central charge $nk$ extending its $V^k$-module structure.
\item Assume that $\mathcal{S}$ satisfies (\ref{Cond1}).
If $V_1=0$ then $\mathfrak{V}(\mathcal{S})_1=0$.
\end{enumerate}
\et
\proof\ The first assertion follows from Proposition \ref{PMax}.

Let $[M]\in R(V)$.
By the definition of the map $w:R(V)\to\{0,1,2\}$ and (b), the lowest weight of $M$ is greater than or equal to $w([M])/2$.
Hence for $[N]\in R(V^k)$ the lowest weight of $N$ is greater than or equal to $w^k([N])/2$.
Thus the latter assertion follows from (\ref{Cond1}) and $V_1=0$.\qe

\br In general, the lowest weight of $M$ for $[M]\in R(V)$ may not be equal to $w([M])/2$.
\er

\br The VOA $\mathfrak{V}(\mathcal{S})$ is rational by \cite{Li}, $C_2$-cofinite by \cite{ABD} and of CFT by the assumption (b).
\er

\br Let $V$ be a simple rational $C_2$-cofinite VOA of CFT type satisfying (a)--(d).
Let $X$ and $N$ be holomorphic VOAs containing $V$ as a full subVOA.
Then by Lemma \ref{LSCE} $X\cong\oplus_{[M]\in\Phi}M$ and $N\cong\oplus_{[M]\in\Psi}M$ as $V$-modules, where $\Phi$ and $\Psi$ are maximal totally singular subspaces of $R(V)$.
By (\ref{EqT}) $\mathfrak{V}(\mathcal{S}(\Phi,\Psi;k))$ is described as follows:
$$\mathfrak{V}(\mathcal{S}(\Phi,\Psi;k))=\bigoplus_{[A(i)]\in\Phi,\ [B]\in\Psi\atop \Pi_{i=1}^k[A(i)]\in\Phi\cap\Psi}\otimes_{i=1}^k (A(i)\times B).$$
\er

\subsection{A description of the moonshine VOA by using $V_{\sqrt2E_8}^+$}
In this subsection, we consider the case where $k=3$ and $V=V_{\sqrt2E_8}^+$, and describe the moonshine VOA as a simple current extension of $V^3$.
We continue the notation of the previous subsection.
For the detail of $V$, see Section 1.4.

Let us describe the moonshine VOA $V^\natural$.
Since $\sqrt2E_8$ has no vectors of norm $2$, we have $V_1=0$.
By Lemma \ref{POVE}, $V$ is a simple rational $C_2$-cofinite VOA of CFT type of central charge $8$ and satisfies (a)--(d) in the previous subsection.
Let $\mathcal{S}$ be a maximal totally singular subspace of $(R(V)^3,q_V^3)$ satisfying (\ref{Cond1}).
Note that such a subspace exists by Proposition \ref{PLM}.
By Theorem \ref{MTV}, $\mathfrak{V}(\mathcal{S})$ is a holomorphic VOA of central charge $24$ with $\mathfrak{V}(\mathcal{S})_1=0$.
Since $V$ is framed, so is $\mathfrak{V}(\mathcal{S})$.
It was shown in \cite{LY} that any holomorphic framed VOA without weight $1$ subspace is isomorphic to the moonshine VOA $V^\natural$.
Hence we obtain the following theorem.

\bt\label{MTVm} Let $V=V_{\sqrt2E_8}^+$ and let $\mathcal{S}$ be a maximal totally singular subspace of $(R(V)^3,q_V^3)$ satisfying (\ref{Cond1}).
Then the VOA $\mathfrak{V}(\mathcal{S})$ is isomorphic to the moonshine VOA $V^\natural$.
\et

\br Let $X$ and $N$ be holomorphic VOAs of central charge $8$ containing $V$ as a full subVOA.
Then by Lemma \ref{LSCE} $X\cong\oplus_{[M]\in\Phi}M$ and $N\cong\oplus_{[M]\in\Psi}M$ as $V$-modules, where $\Phi$ and $\Psi$ are maximal totally singular subspaces of $R(V)$.
Assume that $\Phi\cap\Psi=\{[V]\}$.
Then by Proposition \ref{PLM} and Theorem \ref{MTVm}, 
$$V^\natural\cong\mathfrak{V}(\mathcal{S}(\Phi,\Psi;3))=\bigoplus_{[A(i)]\in\Phi,\ [B]\in\Psi\atop [A(1)]\times[A(2)]\times[A(3)]=[V]}(A(1)\times B)\otimes (A(2)\times B)\otimes (A(3)\times B).$$
as VOAs.
Note that $X$ and $N$ are isomorphic to the lattice VOA $V_{E_8}$ associated to the $E_8$ root lattice by \cite{DM2}.
\er

\br Since the Leech lattice $\Lambda$ contains $(\sqrt2E_8)^3$ as a sublattice, the VOA $V_\Lambda^+$ contains $V^3$ as a full subVOA, and so does $V^\natural$ (\cite{FLM}).
By the decomposition of $V^\natural$ into irreducible $V^3$-modules, one can directly find complementary maximal totally singular subspaces $\Phi_0$ and $\Psi_0$ of $R(V)$ such that $\mathfrak{V}(\mathcal{S}(\Phi_0,\Psi_0;3))\cong V^\natural$ as $V^3$-modules.
Hence Theorem \ref{MTVm} follows from Lemmas \ref{LO} (5), \ref{LSY}, \ref{POVE} (4), Propositions \ref{PSC}, \ref{PVk} (see next subsection) and Theorem \ref{TCh}.
Note that this argument is independent of \cite{LY}.
\er

\bn In \cite{Mi04}, the moonshine VOA was constructed as an extension of $V^3$.
\en

As an application of our description, let us count the dimension of the weight $2$ subspace of $\mathfrak{V}(\mathcal{S})$.
The following lemma is easy (cf.\ \cite{AD}).

\bl\label{LCv} Let $[M]\in R(V)$.
\begin{enumerate}
\item If $w([M])=0$ then the dimension of the weight $2$ subspace of $M$ is $156$.
\item If $w([M])=1$ then the dimension of the weight $1/2$ subspace of $M$ is $1$.
\item If $w([M])=2$ then the dimension of the weight $1$ subspace of $M$ is $8$.
\end{enumerate}
\el

\bp\label{PEqV} Let $V=V_{\sqrt2E_8}^+$ and let $\mathcal{S}$ be a maximal totally singular subspace of $(R(V)^3,q_V^3)$ satisfying (\ref{Cond1}).
Then the dimension of the weight $2$ subspace of $\mathfrak{V}(\mathcal{S})$ is $196884$.
\ep
\proof\ By Proposition \ref{PVk} (see next subsection), $\Aut(V^3)\cong\Aut(R(V)^3,w^3)$.
Hence by Lemma \ref{LSY} and Theorem \ref{TCh}, we may assume that $\mathcal{S}=\mathcal{S}(\Phi,\Psi;3)$.
The weight $2$ subspace of $\mathfrak{V}(\mathcal{S})$ is the direct sum of those of $V^3$ and irreducible $V^3$-modules $M$ with $w^3([M])=4$.
Hence by Lemmas \ref{LN4} and \ref{LCv}, the dimension of the weight $2$ subspace of $\mathfrak{V}(\mathcal{S})$ is $$ (3\times156)+(3\times (2^5-1))\times 1^0\times8^2+(3\times (2^5-1)\times 2^8)\times 1^2\times 8^1=196884.$$\qe 

\subsection{An $E_8$-approach to the automorphism group of $V^\natural$}
In this subsection, as an application of Sections 2.2 and 4.2, we study the automorphism group of the moonshine VOA $V^\natural$.
Note that the results in this section are independent of the fact that the automorphism group of $V^\natural$ is isomorphic to the Monster.
We continue the notation of the previous subsection.

First, we determine the automorphism group of $V^k$.
\bp\label{PVk} Let $V=V_{\sqrt2E_8}^+$.
Then the automorphism group of $V^k$ is isomorphic to $\Aut(V)\wr\Sym_k$.
\ep
\proof\ Clearly, $\Aut(V)\wr\Sym_k$ acts on $V^k$ as an automorphism group.

Let $g\in\Aut(V^k)$.
By Lemma \ref{LemFHL}, $R(V^k)$ is identified with $R(V)^k$.
By Lemma \ref{POVE} (3), $w^k([M])/2$ is the lowest weight of $M$ for any irreducible $V^k$-module $M$.
Hence by Lemma \ref{LAct0} (2) $g\in \Aut(R(V)^k,w^k)$.
By Lemma \ref{POVE} (4) and Proposition \ref{PAutw}, $\Aut(R(V)^k,w^k)\cong O(R(V),q_V)\wr\Sym_k\cong \Aut(V)\wr\Sym_k$.
Hence $g\in \Aut(V)\wr\Sym_k$.
\qe

Next, we consider lifts of symmetries of the quadratic space.
Let $G$ be the automorphism group of $\mathfrak{V}(\mathcal{S})$.
Recall that $\mathfrak{V}(\mathcal{S})$ is a simple current extension of $V^3$ graded by $\mathcal{S}$.
By (\ref{eta}), we obtain a homomorphism of groups
\begin{eqnarray}
\eta:N_{G}(\mathcal{S}^*)\to\{g\in \Aut(V^3)\mid g\circ \mathcal{S}=\mathcal{S}\}.\label{EqVVV}
\end{eqnarray}
By Proposition \ref{PSh2}, $\eta$ is surjective and $\Ker\eta=\mathcal{S}^*\cong 2^{15}$.
We now prove the following theorem without properties of the Monster. 

\bt\label{MT} Let $V=V_{\sqrt2E_8}^+$ and let $\mathcal{S}$ be a maximal totally singular subspace of $(R(V)^3,q_V^3)$ satisfying (\ref{Cond1}).
Let $G$ be the automorphism group of $\mathfrak{V}(\mathcal{S})$ and let $H$ be the stabilizer of $V^3$ in $G$.
\begin{enumerate}[{\rm (1)}]
\item $H=N_{G}(\mathcal{S}^*)$.
\item $G$ is transitive on the set of all full subVOAs of $\mathfrak{V}(\mathcal{S})$ isomorphic to the tensor product of three copies of $V_{\sqrt2E_8}^+$.
\item $H$ is isomorphic to an extension of the stabilizer of $\mathcal{S}$ in $\Aut(R(V)^3,w^3)$ by an elementary abelian $2$-group of order $2^{15}$, and $H$ has the shape $2^{15}.(2^{20}:({\rm SL}_5(2)\times\Sym_3))$.
\end{enumerate}
\et
\proof\ Since the subspace of $\mathfrak{V}(\mathcal{S})$ fixed by $\mathcal{S}^*$ is $V^3$, the normalizer $N_{G}(\mathcal{S}^*)$ is a subgroup of $H$.
Since $\mathfrak{V}(\mathcal{S})$ is a simple current extension of $V^3$, the multiplicity of any irreducible $V^3$-submodule is $1$.
Hence $H$ acts on the set of irreducible $V^3$-submodules as permutations, and $H$ is a subgroup of the normalizer, which proves (1).

Let $X$ be a full subVOA of $\mathfrak{V}(\mathcal{S})$ isomorphic to $V^{3}$.
It suffices to show that there is an automorphism of $\mathfrak{V}(\mathcal{S})$ which sends $X$ to $V^3$.
By Lemma \ref{LSCE}, $\mathfrak{V}(\mathcal{S})$ is a simple current extension of $X$.
Let $\mathcal{T}$ be the set of the isomorphism classes of irreducible $X$-submodules of $\mathfrak{V}(\mathcal{S})$.
Then by Proposition \ref{PMax}, $\mathcal{T}$ is a maximal totally singular subspace of $R(X)(\cong R(V^3))$.
Since the weight $1$ subspace of $\mathfrak{V}(\mathcal{S})$ is trivial, $\mathcal{T}$ satisfies (\ref{Cond1}).
Hence by Theorem \ref{TCh} there exists $g\in\Aut(X)(\cong\Aut(R(V)^3,w^3))$ such that $g\circ\mathcal{T}=\mathcal{S}$.
By Lemma \ref{LSY}, $g\circ \mathfrak{V}(\mathcal{S})$ has a VOA structure isomorphic to $\mathfrak{V}(\mathcal{S})$ and let $\tilde{g}$ be the isomorphism of VOAs from $\mathfrak{V}(\mathcal{S})$ to $g\circ \mathfrak{V}(\mathcal{S})$.
Note that  the set of the isomorphism classes of irreducible submodule of $g\circ \mathfrak{V}(\mathcal{S})$ for $\tilde{g}(X)\cong V^3$ is $g\circ\mathcal{T}=\mathcal{S}$.
Hence by Proposition \ref{PSC}, there is an isomorphism of VOAs $h$ from $g\circ \mathfrak{V}(\mathcal{S})$ to $\mathfrak{V}(\mathcal{S})$ such that $h\circ\tilde{g}(X)=V^3$.
Thus we obtain an automorphism $h\circ \tilde{g}$ of $\mathfrak{V}(\mathcal{S})$ such that $h\circ \tilde{g}(X)=V^3$, which proves (2).

By Proposition \ref{PVk}, $\Aut(R(V)^3,w^3)=\Aut(V^3)\cong \Aut(V)\wr\Sym_3$.
Hence (3) follows from (1), (2), Proposition \ref{PStab} and (\ref{EqVVV}).\qe

Combining Theorems \ref{MTVm} and \ref{MT}, we obtain the following corollary.

\bc Let $V^\natural$ be the moonshine VOA.
\begin{enumerate}[{\rm (1)}] 
\item The automorphism group of $V^\natural$ is transitive on the set of all full subVOAs isomorphic to $(V_{\sqrt2E_8}^+)^{3}$.
\item The stabilizer of a full subVOA of $V^\natural$ isomorphic to $(V_{\sqrt2E_8}^+)^{3}$ in the automorphism group of $V^\natural$ has the shape $2^{15}.(2^{20}:({\rm SL}_5(2)\times\Sym_3))$.
\end{enumerate}
\ec

\br By the order of the Monster (\cite{ATLAS}) and the corollary above, the number of full subVOAs of $V^\natural$ isomorphic to $(V_{\sqrt2E_8}^+)^{3}$ is $$3^{17}\cdot 5^8\cdot 7^5\cdot 11^2\cdot 13^3\cdot 17\cdot 19\cdot 23\cdot 29\cdot 41\cdot 47\cdot 59\cdot 71=391 965 121 389 536 908 413 379 198 941 796 875.$$
\er

\br\label{Rm} Assume that $\mathcal{S}=\mathcal{S}(\Phi,\Psi;3)$.
By Proposition \ref{PStab}, $N_{G}(\mathcal{S}^*)$ normalizes the elementary abelian $2$-subgroup $\Psi_{(123)}^*$ of $\mathcal{S}^*$ of order $2^5$.
\er

\br In \cite{Sh4}, $V^\natural$ was studied as a simple current extension of $V_{\Lambda(5)}^+$ for a certain sublattice $\Lambda(5)$ of the Leech lattice containing $(\sqrt2E_8)^3$.
One can find complementary maximal totally singular subspaces $\tilde\Phi$ and $\tilde\Psi$ of $R(V)$ satisfying $$V_{\Lambda(5)}^+\cong \oplus_{[M]\in{\rm Span}_{\F_2}\{\tilde\Phi_{(12)},\tilde\Phi_{(13)}\}}M\quad {\rm and}\quad \mathfrak{V}(\mathcal{S}(\tilde\Phi,\tilde\Psi;3))\cong V^\natural$$ as $V^3$-modules.
Hence the normalizer of $\tilde\Psi_{(123)}^*(\cong 2^5)$ of shape $2^{5}.2^{16}.2^{2+12}.\Sym_3.{\rm SL}_5(2)$ described in \cite{Sh4} is conjugate to $N_{G}(\mathcal{S}^*)$.
\er

\bn In \cite[Lemma 9.3]{Mi04} a subgroup of shape $2^{7}.2^{20}.2^{12}.({\rm SL}_4(2)\times\Sym_3)$ of the stabilizer of $V^3$ in $\Aut(V^\natural)$ was described as the stabilizer of certain Virasoro frame of $V$ (cf.\ \cite[Section 6.5]{LS2}).
By the shapes, it is a maximal subgroup of $N_{G}(\mathcal{S}^*)$.
\en

\bn According to \cite{ATLAS}, the subgroup of the Monster obtained in Theorem \ref{MT} is a maximal $2$-local subgroup of shape $2^{5+10+20}.({\rm SL}_5(2)\times \Sym_3)$.
\en

\paragraph{\bf Acknowledgements.} The author thanks to Professor Ching Hung Lam for helpful suggestions on this study.
He also thanks to Professor Hiroshi Yamauchi for useful discussion on simple current extensions of VOAs.
Part of the work was done when he was visiting Academia Sinica in 2009 and Imperial College London in 2010.
He thanks the staff for their help.
\begin{small}

\end{small}
\end{document}